\documentclass[a4paper,11pt]{article}
\usepackage{amsmath}

\usepackage{graphicx}

\usepackage{graphics}
\usepackage{amsfonts}

\usepackage{amssymb}
\newcommand{\Rs}{\mathbb{R}}
\newcommand{\Ns}{\mathbb{N}}
\newcommand{\Zs}{\mathbb{Z}}

\newcommand{\lp}{\left(}

\newcommand{\rp}{\right)}

\newcommand{\noi}{\noindent}
\newcommand{\beq}{\begin{equation}}
\newcommand{\eeq}{\end{equation}}

\newtheorem{THE}{Theorem}[section]

\newtheorem{LE}[THE]{Lemma}
\newtheorem{PROP}[THE]{Proposition}
\newtheorem{DEF}[THE]{Definition}
\newtheorem{COR}[THE]{Corollary}

\newtheorem{RMK}[THE]{Remark}

\newtheorem{PROP-DEF}[THE]{Proposition-Definition}

\begin{document}

\title{Morse Novikov Theory and Cohomology with Forward Supports}
\author{F. Reese Harvey (harvey@math.rice.edu) \and Giulio Minervini (minervini@dm.uniba.it) }
\maketitle

\begin{abstract}

We present a new approach to Morse and Novikov theories, based on the deRham Federer theory of currents, using the finite volume flow technique of Harvey and Lawson [HL]. In the Morse case, we construct a noncompact analogue of the Morse complex, relating a Morse function to the \emph{cohomology with compact forward supports} of the manifold. This complex is then used in Novikov theory, to obtain a geometric realization of the Novikov Complex as a complex of currents and a new characterization of Novikov Homology as cohomology with compact forward supports. Two natural ``backward-forward'' dualities are also established: a \emph{Lambda duality} over the Novikov Ring and a \emph{Topological Vector Space duality} over the reals.
\end{abstract}

\section*{Introduction}

A new approach to Morse theory on compact manifolds has recently been introduced by Harvey and Lawson in [HL]. The goal of this paper is to construct a similar approach to Novikov theory. With this aim, it is fundamental to first develop a version of Morse theory for nonproper Morse functions, which is interesting in itself. 

The present paper thus consist of two parts. To better describe its contents, it is worthwhile to recall some facts about Morse and Novikov theories. It is not our intention to give a complete or detailed history; we instead confine the comments to what is pertinent to this paper.\\

\subsection*{Part 1 : Non-proper Morse Theory} We recall that a function $f : Y \rightarrow \Rs$ is \textbf{Morse} if all its critical points $x$ are not degenerate. The \textbf{index} $\# (x)$ of the critical point $x$ is the dimension of a maximal subspace where the Hessian of $f$ is negative definite. If $\phi_t$ is a gradient flow for $-f$, each critical point $x$ is a fixed point of the flow, and one defines the \textbf{Stable }and \textbf{Unstable manifold} at $x$ as
\[S_x = \left\{ p\in Y \, |\, \lim_{t\rightarrow + \infty} \phi_t (p) = x  \right\} \ \ \ \ \  \mbox{ and } \ \ \ \ \  U_x = \left\{ p\in Y \, |\, \lim_{t\rightarrow - \infty} \phi_t (p) = x \right\} \]
The Stable Manifold theorem states that $S_x$ is an injectively immersed copy of $\Rs^k$, where $k = n-\#(x)$ and $n=$dim$(Y)$. Similarly $U_x$ is a submanifold of dimension $\#(x)$, which is transversal to $S_x$.\\

The idea that the critical points of a Morse function $f:Y\rightarrow \Rs$ describe the topology of the manifold $Y$, when $Y$ is compact, is very old. However, some formulations of this idea  are quite recent, and still in progress. As explained in [Bo], it was known after Thom [T] that the stable manifolds of the critical points of $f$ provide a cellular partition of $Y$, but the singularities on their boundary were not well understood, and the usual tools of Morse theory have been the change of homotopy of prelevel sets and the Morse inequalities.  

After the papers of Witten [W] and Floer [Fl], it has become customary to encode the global information of the Morse function in the so called \textbf{Morse Complex} $(C_{\mathcal{M}}^*, \delta)$, defined using a Smale gradient for the function $f$. The group $C^k_{\mathcal{M}}$ is the free abelian group generated by the critical points of index $k$ of $f$, which are finite since the manifold $Y$ is compact. The differential $\delta $ is defined by ``counting with orientation'' the flow lines connecting critical points whose index differs by 1, which are finite because the flow is Smale. For the rule for counting flow lines cf. the given references or Remark \ref{RMK:Counting Rule} below.\\ 

\textbf{Fundamental Theorem of compact Morse theory} \newline
\emph{The Morse complex $C_{\mathcal{M}}^*$ is indeed a complex (i.e. $\delta^2 =0$) and its cohomology $H^k (C^*_{\mathcal{M}} )$ is isomorphic to the standard cohomology $H^k (Y, \Zs)$.}\\ 

The classical Morse inequalities are a purely algebraic corollary, but the previous theorem gives substantial additional information because of the explicit knowledge of the differential $\delta$. The Morse complex is thus a convenient tool to state Morse theory, but remains a formal complex: the ``cell'' associated to each critical point is disguised.

The first geometric realization of this cell and of the Morse complex is due to Laudenbach, [La]. Under a technical ``tameness'' assumption on the flow (to be linear near the fixed points, with eigenvalues $=\pm 1$) he resolves the singularities of the stable manifold $S_x$. This allows one to introduce the current of integration over $S_x$ and to prove that the rule for counting flow lines in the Morse complex is related to computing its current boundary (cf. [La] or Corollary \ref{COR: boundary of stable2} below). Replacing $x$ by $S_x$, the Morse complex is thus realized as the subcomplex $\mathcal{S}^*$ of the complex of currents, (finitely) generated by (the currents defined by) the stable manifolds. We will refer to the complex $\mathcal{S}^*$ as the \textbf{current Morse complex}.
     
Harvey and Lawson [HL] go further, introducing the technique of ``finite volume flows''. Again assuming the tameness condition above, they show that the limits of the pullbacks under the gradient flow $\phi_t$ of any smooth form $\alpha$:
\begin{equation}\label{eq:intro1}
\lim\limits_{t\rightarrow+\infty}\phi_{t}^{\ast}\left(  \alpha\right) = \mathbf{P}\left(  \alpha\right) = \sum\limits_{x\in Cr\left(  f\right)  } \left(  \int_{U_{x}}\alpha\right)  S_{x} 
\end{equation}
converge (as currents) to chains in $\mathcal{S}^*$. They construct in addition an explicit operator $\mathbf{T}$, providing the following chain homotopy between the limit operator $\mathbf{P}$ and the identity on smooth forms $\alpha$:   
\begin{equation}\label{eq:intro2}
d \mathbf{T} \alpha +\mathbf{T} d \alpha =\mathbf{\alpha }-\mathbf{P}\alpha   \ \ \ \ \ \mbox{where} \ \ \ \ \  \mathbf{T}(\alpha)  = \int_0^{+ \infty} \phi^*_t (\alpha)
\end{equation}
This yields a deformation of the deRham complex of forms onto the Current Morse Complex $\mathcal{S}^*$; replacing smooth forms by chains, the result extends to integer coefficients, providing a new, more transparent proof of the Fundamental Theorem above.\\ 

In this paper, we study the main question of Morse theory when the manifold is not necessary compact and the function not necessarily proper.\\ 

\noi \textbf{Q }. How does a Morse function relate to the topology of the manifold?\\ 

We do assume a \emph{Weakly Proper} condition, motivated by the examples arising from Novikov theory, to prohibit the flow from concentrating an infinite amount of mass in small regions, cf. Definition \ref{DEF:Weakly-Proper}.

Under this assumption, we show that for any test form $\alpha$, the limit 
\[\lim\limits_{t\rightarrow+\infty}\phi_{t}^{\ast}\left(  \alpha\right) = \mathbf{P}\left(  \alpha\right) = \sum\limits_{x\in Cr\left(  f\right)  } \left(  \int_{U_{x}}\alpha\right)  S_{x}\]
still converges, as is the case for compact manifolds. In spite of the compact support of $\alpha$, the limit may not be compactly supported, since the stable manifolds need not be relatively compact. This unfortunately makes the homotopy between $\alpha $ and $\mathbf{P} (\alpha)$ not directly fruitful.
 
Nevertheless, observing that the flow moves in a preferred direction motivates the introduction of the family of ``compact forward'' subsets and the definition of the ``$\mathcal{S}^*$ complex with compact forward supports'', denoted $\mathcal{S}_{c \uparrow}^*$. This complex is an appropriate analogous of the Current Morse Complex: it is made up of currents which are sums of stable manifolds and have compact forward support. In particular $\mathcal{S}_{c \uparrow}^*$ contains the limits of test forms $\mathbf{P}(\alpha)$.\\

The main Theorems \ref{THE:Morse-reals} and \ref{THE:Morse-integers} of the first part describe our solution to the question \textbf{Q}. They state, in particular, that the cohomology of the complex $\mathcal{S}_{c \uparrow}^*$ is isomorphic to the cohomology of $Y$ with supports in the family of compact forward sets (with either real or integer coefficients).\\

Although it is not widely used, the ``cohomology of $Y$ with support in a family of subsets of $Y$'' is a standard tool in sheaf cohomology (cf. [G]). The family of compact forward sets is sufficiently nice so that, for example, choosing real coefficients, the groups $H_{c \uparrow }^k (Y , \Rs )$ are isomorphic to the deRham groups defined using smooth forms with compact forward supports.\\

The exposition in this first part is brief and it largely follows the lines in [HL]. In particular Theorem \ref{THE:M} is given without proof, since the proof is not fundamental for the understanding of the rest of the paper and will appear in [M2] in a greater generality, not needed for the applications to Novikov theory. A complete exposition can be found in the second author's PhD thesis [M].\\

\subsection*{Part 2 : Novikov theory} Novikov theory appeared in the early '80 as a generalization of Morse theory to multivalued functions and closed 1-forms. In [N], Novikov considers the problem of finding relations between their critical points and the topology of the manifold. Next is a brief outline of his construction.

Given a closed 1-form $\omega$ on a compact manifold $X$ (for a function $g: X \rightarrow S^1$ take $\omega =dg$), one utilizes a covering 
\[
\pi: Y \longrightarrow X
\]
where the pullback $\pi ^* (\omega ) = df$ is an exact form; a genuine Morse function $f$ thus arises on the noncompact manifold $Y$. 

A complex, say $(C^k_{\mathcal{N}}, \delta)$, now called the \textbf{Novikov complex}, is then defined. It is similar, ante litteram, to the current Morse complex and at the same time a generalization of it. The generators of $C^k_{\mathcal{N}}$ are in 1-1 correspondence with the critical points of the form $\omega$ on $X$. In fact, Novikov describes the generator associated to a critical point $x$ for $\omega$ as the ``infinite ascending cell'' (the stable manifold, in our terminology) of any point $y$ which lifts $x$ to $Y$, and hence is critical for $f$.
The boundary $\delta$ can be defined by counting flow lines, as in the Morse complex. However, the computation involves infinite contributions (but locally finite), because there are infinitely many critical points, cf. e.g. [Sc2], the definition in [N] being somewhat elusive.         

Using the ingenious trick of choosing coefficients in a ring of formal series (the \emph{Novikov Ring}), Novikov is able to make the complex $C^k_{\mathcal{N}}$ and its cohomology groups $H^k(C^*_{\mathcal{N}})$ (called \textbf{Novikov homology}) into finitely generated modules, and to obtained the \emph{Novikov inequalities}, relating the ranks of the Novikov homology modules to the number of critical points of $\omega$.\\

The Novikov complex has not been very well understood, mainly because the tools needed for a coherent analysis of the theory were not completely developed. The problem of the singularities of the cell (already encountered in compact Morse theory), added to the problems of infinite support and of taking infinite summations. 

Many authors thus felt the need of different approaches to Novikov theory, interpreting the Novikov complex within standard algebraic topology.     

The first and still most cited such result is due to Pajitnov, [P3]. He proves that the Novikov complex is chain homotopic to the complex $\Lambda \otimes_{\Zs [\pi ]} C^ * (Y)$, where $\Lambda$ is the Novikov ring, $\Zs [\pi ]$ is the group ring of the covering $ \pi :Y\rightarrow X$ and $C^* (Y)$ is any standard cochain complex for $Y$, invariant under the deck translations (e.g. the simplicial complex associated to some triangulation of $Y$ obtained by lifting a triangulation of $X$).\\

In this paper the machinery of Part 1 is used to obtain a quick development of Novikov theory, very close in spirit to Novikov's original formulation. In particular we give a geometric realization of the Novikov complex, a new characterization of Novikov homology, and a new, short proof of Pajitnov's theorem.\\

In the geometric framework above, we first construct a Weakly Proper gradient flow for $f$ on $Y$. The Novikov ring $\Lambda$ can be regarded as acting on subsets of $Y$, and this action commutes with the flow and preserves the family of compact forward sets. The complex $\mathcal{S}^*_{c \uparrow}$ thus naturally becomes a module over the Novikov ring $\Lambda$. 

As a complex of $\Lambda$-modules, the complex $ \mathcal{S}_{c\uparrow}^{\ast} $ is isomorphic to the Novikov complex, see Proposition \ref{PROP:current Novikov complex}. The Novikov complex is thus geometrically realized as a subcomplex of currents, just as the current Morse complex realizes the Morse complex.\\

The main Theorems \ref{THE:main} and \ref{THE:main-general} follow directly from Part 1, by taking the cohomology of $ \mathcal{S}_{c\uparrow}^{\ast} $. They provide the characterization of Novikov homology as cohomology with compact forward supports.\\

We then observe that the family of compact forward set does not vary for bounded perturbations of the defining function. This allows us to modify the function $f$ by a suitable (large but bounded) perturbation, in order to get a simpler Novikov complex, following an idea of Latour. After this modification, it is straightforward to obtain Pajitnov's theorem with our method, see Theorem \ref{THE:pajitnov}.\\     

The present current approach to Novikov theory has several advantages. First, the presentation is rapid and clean, since we are able to treat the stable manifolds as geometric chains, without the necessity of introducing auxiliary triangulations etc. As a good test, note that the exposition in Section 5 of the general case of Novikov theory for 1-forms does not essentially differs from the more elementary case of circle valued functions in Section 4.

Second, for the first time we clearly distinguish the Morse-theory ingredients from the algebra involved in the construction: our definition of the Novikov complex requires no algebra and has an intrinsic geometric and topological content. The Novikov ring is only used to obtain finitely generated invariants.
 
Third, the strength of the finite volume approach to Morse theory extends to Novikov theory: the explicit (Morse chain) homotopy between the geometric complexes of smooth forms or chains and the Novikov complex is made available.\\

An analogue for Novikov theory of the Fundamental Theorem of Morse theory involves a choice, resulting in several different statements. Pajitnov's theorem is one, and our main Theorems \ref{THE:main} and \ref{THE:main-general} provide a second. Still another possibility is contained in Farber-Ranicki [FR]. It is important to note that for any one of these results the ensuing inequalities are an algebraic corollary, as is the case for Morse theory. The Novikov inequalities were a tough subject in early papers (the first complete proofs are probably due to Farber [F] and Pajitnov [P2], who notably do not seem to utilize the Novikov complex) but they are now simplified (cf. Remark \ref{RMK:Lambda-inequalities} or the discussion in [R], p.12). \\

In the last section, we discuss two ``backward-forward'' pairings, leading to dualities for Novikov homology: a topological vector space duality over the reals and a ``Lambda'' duality over the Novikov Ring $\Lambda$. The former is new, whereas a weaker version of the latter is due to Pajitnov [P2], cf. also [Sc2].  

An algebraic appendix dealing with the Novikov ring and related rings of formal series is added for completeness.

\section{Weakly-Proper Smale flows}

Suppose $f:Y\rightarrow\mathbb{R}$ is a Morse function on an oriented 
(not necessarily compact) Riemannian manifold and assume that the vector field $V=-grad(f)$ is complete, insuring a flow $\phi= (\phi_t )_{t \in \Rs}$ on $Y$. Recall the transversality hypothesis of Smale:

\begin{DEF}The flow is \textbf{Smale} if for any two critical points $x,y\in
Cr\left(  f\right)  $, the stable manifold $S_{x}$ and the
unstable manifold $U_{y}$ intersect transversally.
\end{DEF}

To overcome lack of compactness of $Y$, we will also assume the following condition (see [M2] for a generalization to nongradient flows):

\begin{DEF}\label{DEF:Weakly-Proper} The flow is \textbf{Weakly Proper} if the intersection of each broken flow line with each slab $f^{-1}\left(  \left[  a,b\right]  \right)  $ is compact.
\end{DEF}

Of course any gradient flow of a proper function is weakly proper, but this case is not general enough for our purposes. \\

Denote by $\Phi$ the \textbf{total graph map}
\[
\begin{array}{ccc}
\Phi :\mathbb{R} \times Y  \rightarrow   Y\times Y  &\text{ sending }& (t,p)  \mapsto  (\phi_{t}\left(  p\right)  ,p)
\end{array}
\]

This map is regular near points $(t,p)$ such that $p \notin Cr(f)$. Using $\Phi$, we can define the following pushforwards:
\begin{equation*}
P_{t}=\Phi_{\ast}(t\times Y)=\left[  graph(\phi_{t})\right]\text{ \ \ \ and\ \ \ \ } T_{t} = - \Phi_{\ast}(\left[  0,t\right]  \times Y) 
\end{equation*}
These define currents of integration in $Y\times Y$ of dimension respectively $n= \text{dim}Y$ and $n+1$; in particular $graph(\phi_{0})$ is just the diagonal $\Delta$. Since taking boundary commutes with the current pushforward:
\begin{equation*}
\partial T_{t}= \Delta - P_{t} \tag{$\textbf{Stokes}$\ $\textbf{Theorem}$}\end{equation*}

\begin{RMK}\label{RMK:degree-dim} Here and in the sequel $\partial$ denotes the current boundary, i.e. the adjoint of the exterior derivative (differential) on forms. It has the disadvantage that $\partial \alpha = (-1)^{k\!+ \! 1} d \alpha$ for smooth $k$-forms. Later we will also use the ``differential'' $d$ on currents, defined by $d R =(-1)^{\text{deg} R \! + \! 1} \partial R$ , which extends the exterior differential on forms.
\end{RMK}  

Consider the 1-1 immersed submanifold $T= \Phi \lp (0,\! +\infty )  \times \lp Y \! - C \! r\left(\!  f \! \right) \rp \rp$.

\begin{DEF}
If the submanifold $T$ has locally finite volume in $Y\times Y$ then the gradient flow $\phi$ is called a \textbf{Finite Volume Flow }(this is actually a property of the associated singular foliation).
\end{DEF}

A more general concept was introduced in [HL], without assuming the flow to be a gradient. In this case, the volume of $T$ has to be computed with ``multiplicities'' in order to take into account periodic behavior.

Local bounds on the volume of $T$ can be obtained by looking at the closure $\overline{T}$, which has a natural structure of a stratified set. In [M] the following concept was introduced:

\begin{DEF} A stratified space $S\subset Y$ is called \textbf{Horned} provide it is AB Whitney regular and there exists a locally finite covering of \emph{desingularizations}, each consisting of a compact manifold with corners $M$ of the same dimension as $S$ and a stratified map $\pi :M \rightarrow S$ whose restrictions to each stratum are submersions.
\end{DEF}
 
Any horned stratified space $S$ has locally finite volume, since the desingularizations provide local bounds for the volume. Therefore, if $S$ is oriented, it defines a current of integration, denoted by $[S]$ (or simply $S$ if no confusion may arise). 

If $U,S \subset Y$ are horned, oriented and transversal, i.e. all possible intersection of strata are transversal, then $U \cap S$ is horned and canonically oriented (since $Y$ is oriented). With this convention, it is always true that $[U \cap S] = [U] \wedge [S]$.\\      

In this paper we utilize one of the main results in [M].

\begin{THE}\label{THE:M} Suppose the gradient flow is Weakly Proper and Smale.

  The closure of the stable (resp. unstable) manifold of a critical point $x$ is a horned stratified space, whose topological boundary is contained in the union of the stable (resp. unstable) manifolds of critical points $y$, with $\#(y) \! > \! \#(x) $ and $f(y)\! > \! f(x)$ (resp. $\#(y) \! < \! \#(x) $ and $f(y)\! < \! f(x)$).  

The closure of the submanifold $T \subset Y\times Y $ is a horned stratified space and the flow is thus finite volume. The topological boundary of $T$ is made up of the diagonal $\Delta$ and the union of the closures of the products $U_x \times S_x$, over all the critical points $x$.
\end{THE}

The stable manifolds or the submanifold $T$ are special cases of ``shadows''. The shadow of a submanifold $K$ under a flow is the union of the trajectories starting on $K$: its singularities arise from the singular points of the flow. Under the hypothesis of the theorem, it can be proved that the shadow of any submanifold $K$ contained in a level set and transversal to all unstable manifolds is horned stratified. 
    
On compact manifolds, similar results were proved in [La] and [HL] under a special condition on the fixed points of the flow (to be linear with eigenvalues $=\pm 1$). The additional hypothesis is not assumed in [AB], where Morse-Smale (and other) gradients are discussed from a slightly different viewpoint, along with an important ODE boundary value technique. This technique, suggested by L. Shilnikov in the sixties (cf. the monography [S]), provides control of the otherwise elusive asymptotics of a shadow near a critical point. The weakly proper condition avoids concentrations, so one can proceed iteratively, in desingularizing shadows, by climbing the critical points encountered. See [M] or [M2] for a complete proof and generalizations.\\

The orientation on $T$ is the opposite of that induced by $\Phi$, whereas we choose arbitrary orientations for each stable manifold $S_x$ and orient the corresponding unstable manifold $U_x$ so that $<\! U_x \! > \! + \! <\! S_x \! > = <\! Y\! \! >$. 

Since the volume of $T-T_{t}$ locally decreases to zero, the limit of currents
\[
\lim_{t\rightarrow+\infty}T_{t} = T
\]
holds (in the mass norm) and the family $P_{t} = \Delta -\partial T_{t}$ necessarily also converges.  A simple geometric description can be established for the limit $P$, generalizing Theorem 3.3 in [HL], namely:

\begin{COR}\label{COR:FME}
The currents $P_{t} = [graph ( \phi_t ) ]$ on $Y \times Y$ converge to: 
\begin{equation}\label{current P}
\lim\limits_{t\rightarrow
+\infty}P_{t} \overset{def}{=} P =
\sum_{x\in Cr\left(  f\right)  }
\left[  U_{x}\times S_{x}\right] 
\end{equation}
and the following ``\textbf{Fundamental Morse Equation}'' holds: 
\begin{equation*}
\partial T=\Delta-P \tag{$\textbf{FME}$}
\end{equation*}
\end{COR}

The use of the current boundary $\partial$ in (FME), instead of the differential $d$, is coherent with the rules of the kernel calculus, adopted in the next section.

A second consequence of Theorem \ref{THE:M} is the extension of a basic result, first pointed out by [W] and stated, in the following geometric form, in [HL], Proposition 4.5 (see also [La]). 

\begin{COR}\label{COR: boundary of stable2}
Let $x$ be a critical point and $S_x$ its stable manifold. Then the boundary of the current of integration over $S_x$ is:
\begin{equation}\label{eq:boundary of stable2}
d S_x = \sum\limits_{\# y = 1 + \# x}  n_{x}^{y} \, S_{y}
\end{equation}
Each coefficient $n_{x}^{y} \in \Zs$ is determined by ``counting with orientation'' the flow lines connecting $x$ to $y$.  
\end{COR}

The Smale hypothesis is fundamental for formula (\ref{eq:boundary of stable2}) and proves that the flow lines connecting $x$ to a fixed $y$ with index $1 + \# (x)$ are ``discrete''; the Weakly Proper condition then insures that they are actually finite and that the sum in (\ref{eq:boundary of stable2}) is locally finite.

Intuitively, for any $x$ and $y$ as above, the stable manifold $S_x$ accumulates along $S_y$ near $y$ like (finite) pages in a book, having $S_y$ as common edge, and each page $F_j$ satisfies $dF_j = \pm S_y$. Intersecting $S_x$ with $U_y$ selects one flow line $\gamma_j$, ending at $y$, in each page $F_j$ and $\gamma_j$ is to be counted positively or negatively according to the $\pm$ sign above.\\ 

The above argument is made rigorous in [HL] and [La] using the fact that the singularity of $S_x$ near $y$ is that of a manifold with boundary. The same proof establish formula (\ref{eq:boundary of stable2}) in the present case, because the singularity is ``horned''. Since $Y$ is oriented, the rule for counting flow lines can be determined using the Leibniz rule as follows.\\

Let the orientation on $U_y \wedge  S_x = \sum \epsilon_j \gamma_j$ define the $\pm$ signs $\epsilon_j$, and put $k  = \mbox{deg}(U_y)$. Near $y$, where $d (U_y) =  0$ and $d\gamma_j =  (-\! 1)^{n \! + \! 1 }[y]$ (cf. Remark \ref{RMK:degree-dim}), we can compute:
\begin{equation*}
n_{x}^{y} \, [y] = U_y \wedge d (S_x) = (-1)^{k} d (U_y \wedge S_x) =  (-1)^{k} \sum (-\! 1)^{n \! + \! 1} \epsilon_j   [y] 
\end{equation*}

Since $k = \mbox{deg}U_y = n \! - \! \#(x) \! - \! 1 $, this proves: 
\beq
n_{x}^{y} =  (-1)^{\#(x)} \sum \epsilon_j
\eeq

\begin{RMK}[Counting Rule]\label{RMK:Counting Rule}The flow line $\gamma_j$ ``counts'' positively, either when its orientation agrees with the orientation of $U_y \cap S_x$ and $\#(x)$ is even, or when the orientations disagree and $\#(x)$ is odd.
\end{RMK}

\section{Forward Supports: the Current Morse Complex}\label{SEC:Forward}

Recall the kernel calculus developed by Harvey and Polking in [HP]. Any ``kernel'' current $R$ on $Y\times Y$ defines, by partial integration, an operator $\mathbf{R} : \mathcal{E}^*_{cpt} (Y) \rightarrow \mathcal{D'}^{*} (Y) $ mapping test forms on $Y$ to currents on $Y$ (the notation $\mathcal{D'}^{k} (Y)$ is standard, but we use $\mathcal{E}^{k}_{cpt} (Y)$ instead of the usual $\mathcal{D}^{k} (Y)$ because of our use of sheaf theory: $\mathcal{D}^k$ is ``not a sheaf'').\\

More precisely, if $\alpha$ is any test form on $Y$ and $\pi_1 ,\pi_2$ the two projections $Y \times Y \rightarrow Y$, one defines $\mathbf{R} (\alpha) = \pi_{2*} \lp R \wedge \pi_1^* (\alpha) \rp  $. In other words:
\[
\lp \mathbf{R}(\alpha)\rp (\beta) = R\lp \pi_1^* (\alpha) \wedge \pi_2^* (\beta)\rp 
\]

\noi for any test form $\beta$ on $Y$, of appropriate degree. 

Under this identification of kernels $R$ and operators $\mathbf{R}$, the diagonal $\Delta$ determines the identity operator $\mathbf{I}$ whereas the graphs $P_{t}=\left[  graph\left(  \phi_{t}\right)  \right]  $ determine the pullback operators $\mathbf{P}_t = \phi^*_t$. The operator $\mathbf{T}: \mathcal{E}_{cpt}^{k}\left(  Y\right) \rightarrow \mathcal{D'}^{k-1}\left(  Y\right)$ corresponding to the kernel current $T$ of the previous section is:
\begin{equation}\label{eq:T}
\mathbf{T}(\alpha) = \int_0^{+ \infty} \phi^*_t (\alpha) 
\end{equation}
whereas the operator $\mathbf{P}: \mathcal{E}_{cpt}^{k}\left(  Y\right) \rightarrow \mathcal{D'}^{k}\left(  Y\right)$, corresponding to the kernel current $P$ is given by:
\begin{equation}\label{eq:P} 
\lim\limits_{t\rightarrow+\infty}\phi_{t}^{\ast}\left(  \alpha\right) = \mathbf{P}\left(  \alpha\right)  = \sum\limits_{x\in Cr\left(  f\right)  } \left(  \int_{U_{x}}\alpha\right)  S_{x} 
\end{equation}
Here $\alpha$ is any test form on $Y$ and convergence is as currents. Finally, the operator $d\circ \mathbf{T}+\mathbf{T}\circ d$ corresponds to the kernel $\partial T$. Using these correspondences, Corollary \ref{COR:FME} determines the following equation of operators:
\begin{equation*}
d\circ \mathbf{T}+\mathbf{T}\circ d=\mathbf{I}-\mathbf{P} \tag{$\textbf{MCH}$}\label{eq:MCH}
\end{equation*}
which will be referred to as the \textbf{Morse Chain Homotopy}.\\ 

It is important to note that the limit $\lim \phi_{t}^{\ast}\left(  \alpha\right) = \mathbf{P}\left(  \alpha\right)$ of a form with compact support is not smooth nor compactly supported, but on the contrary is the sum of stable manifolds. As such, it has ``compact/forward'' support. 

\begin{DEF}\label{DEF:forward} A closed set $A\subset X$ is a \textbf{compact/forward} set (abbreviated c/f set) with respect to the function $f$ if both

\begin{itemize}
\! \item $A\cap f^{-1}\left(  \left[  b,c\right]  \right)  $ \emph{is compact
for any }$b\leq c\in\mathbb{R}$\emph{ (i.e. }$A$ \emph{is \textbf{slab
compact})}

\! \item $A\subset f^{-1}\left(  [a,+\infty)\right)  $\emph{ for some constant
}$a\in\mathbb{R}$\emph{ (i.e. }$A$\emph{ is \textbf{forward}).}
\end{itemize}
\end{DEF}

One can define \textbf{backward} and \textbf{compact/backward} sets in a similar manner. Note that a closed set is compact/forward if and only if it has compact intersection with all the backward sets, and so on for the other cases. The subscript $\uparrow$ will denote the family of forward sets, while $c\uparrow$ or $c/f$ will denote the compact/forward family. For example, $\mathcal{E}_{c\uparrow}^{\ast}\left(  Y\right) =\Gamma_{c\uparrow}(Y,\mathcal{E}^{\ast})$ stands for the space of smooth forms with c/f support.

\begin{RMK}\label{RMK:forward} In the rest of the section we will only consider compact/forward sets (since primarily these will be used in the second part), but everything holds by replacing compact/forward with forward with obvious modifications.
\end{RMK}

\begin{LE}
The operators $\mathbf{T}$ and $\mathbf{P}$ extend to mappings from $\mathcal{E}_{c\uparrow}^{\ast}\left(  Y\right)  $ to  $\mathcal{D'}_{c\uparrow}^{\ast}\left(  Y\right)$ and the Morse Chain Homotopy equation continues to hold. In particular $\mathbf{T}$ is a chain homotopy between the operators
\[
\mathcal{E}_{c\uparrow}^{\ast}\left(  Y\right)  \overset{\mathbf{P}}{\longrightarrow} \mathcal{D}_{c\uparrow}^{\prime\ast}\left(  Y\right)
\, \, \, \text{and} \, \, \,  \,
\mathcal{E}_{c\uparrow}^{\ast}\left(  Y\right)  \overset{\mathbf{I}}{\longrightarrow} \mathcal{D}_{c\uparrow}^{\prime\ast}\left(  Y\right)
\]
\end{LE}

We next introduce the $\mathcal{S}$-complex of currents:

\begin{DEF}The (\textbf{compact/forward}) \textbf{$\mathcal{S}$-complex over $\Zs$}, denoted by $ {}_{\mathbb{Z}}\mathcal{S}_{c\uparrow}^{\ast
}\left(  f \right)  $ is the subcomplex of $\mathcal{D}_{c\uparrow
}^{\prime\ast}\left(  X\right)  $ (the complex of currents with c/f
support) consisting of those currents of the form 
\[
\sum_{x\in F} a_{x}\left[  S_{x}\right] \text{ \emph{ where F  is a c/f set of
critical points and }} a_{x}\in\mathbb{Z}
\]
The boundary $d:{}_{\mathbb{Z}}\mathcal{S}_{c\uparrow}^{\ast}\left(
f \right)  \rightarrow {}_{\mathbb{Z}}\mathcal{S}_{c\uparrow}^{\ast
}\left(  f \right)  $ is the differential on currents.

Similarly we define the $\mathcal{S}$-complex with real coefficients $_{\mathbb{R}}\mathcal{S}_{c\uparrow}^{\ast}\left( f \right)  $ .
\end{DEF}

Each element of $\emph{{}}_{\mathbb{Z}}\mathcal{S}_{\uparrow}^{\ast}\left(  f\right)  $ is thus the sum of a locally finite family of currents defined by stable manifolds. Note that the $\mathcal{S}$-complex depends on the function $f$, which defines the c/f supports, and on the Riemannian metric on $Y$, which determines the gradient field and hence the stable manifolds. The symbol $f $ will be deleted from $\emph{{}}_{\mathbb{Z}}\mathcal{S}_{c\uparrow}^{\ast}(f)$ when clear form the context; the Riemannian metric is instead always understood.\\

The compact/forward sets trivially form a paracompactifying family for $Y$, and one can consider sheaf cohomology with compact/forward supports. This support family is nice enough so that $H_{c\uparrow}^{\ast}\left(  Y,\mathbb{R}\right)$ can be computed using either the complex $\mathcal{D}_{c\uparrow}^{\prime\ast}\left(  Y\right)  $ of currents with c/f support, or the complex $\mathcal{E}_{c\uparrow}^{\ast}\left(  Y\right)  $ of smooth forms with c/f support (see [G], Theorems 3.5.1 and 4.7.1). For a discussion about computing cohomology using various sheaves of currents, see the appendix in [HZ].

Now we can state the main results of the first part, Theorems \ref{THE:Morse-reals} and \ref{THE:Morse-integers}. They relate the cohomology of the compact/forward $\mathcal{S}$ complex to the topology of the manifold, and provide an answer to the question \textbf{Q} in the introduction.

\begin{THE}\label{THE:Morse-reals} Suppose the gradient flow is Weakly Proper and Smale. Then the maps 
\[
\mathbf{P}:\mathcal{E}_{c\uparrow}^{\ast}\left(  Y\right)  \mathcal{\longrightarrow
}_{\mathbb{R}}\mathcal{S}_{c\uparrow}^{\ast} \mbox{     and     }\mathbf{I}:\  _{\mathbb{R}}\mathcal{S}_{c\uparrow}^{\ast}
\hookrightarrow\mathcal{D}_{c\uparrow}^{\prime\ast}\left(  Y\right)
\]
induces (algebraic) isomorphisms in cohomology:
\begin{equation}
H_{c\uparrow}^{k}\left(  Y,\mathbb{R}\right)  \approx H^{k}\left(
_{\mathbb{R}}\mathcal{S}_{c\uparrow}^{\ast}  \right) \nonumber
\end{equation}
\end{THE}

\noi \textbf{Proof}. Since $\mathbf{I} \circ \mathbf{P}$ is homotopic to the inclusion $\mathcal{E}_{c\uparrow}^{\ast}\left(  Y\right)  \hookrightarrow \mathcal{D}_{c\uparrow}^{\prime\ast}\left(  Y\right)$ and the latter induces an isomorphism in cohomology (by deRham's theorem), it suffices to show that the map induced in cohomology by $\mathbf{P}$ is surjective. That is, if $S\in {}_{\mathbb{R}}\mathcal{S}_{c\uparrow}^{\ast}$ is closed, then there is a closed form $\varphi \in \mathcal{E}_{c\uparrow}^{\ast}\left(  X\right)$ such that $P(\varphi)=S$.
This fact is a little technical, though not very difficult and the proof is a refinement of an argument in [HL], p.12.
Suppose $S\in {}_{\mathbb{R}}\mathcal{S}_{c\uparrow}^{k} $ is closed and $f(S) \geq 0$. The idea is to consider slabs $f^{-1}([0,n])$ and find a sequence of forms $\varphi_{n}$  with compact/forward support on $X$ such that:\\

\noindent- $\varphi_{n}$  is closed \newline
\noindent- spt$(\varphi_{n})\subset f^{-1}([n,+\infty[)$ and it is compact/forward \newline
\noindent- $\mathbf{P}\lp\sum_{i=1}^{n}(\varphi_{i})\rp - S$ vanishes on $f^{-1}([0,n])$ \\

See [M2] for the details of this construction. Then $\varphi = \sum \varphi_{n}$ defines a closed form with c/f support and $P\lp\varphi\rp= S$ $\Box$\\

Theorem \ref{THE:Morse-reals} can be strengthened with integer coefficients. One can use the $\mathcal{S}$-complex over $\mathbb{Z}$ to compute $H_{c\uparrow}^{\ast}\left(  Y,\mathbb{Z}\right)$ as follows. 
After deRham ([D]), a \textbf{local chain current} is a current that can be locally described as a finite sum of (currents defined via pushforward by)
smooth simplexes. Let's denote by $\mathcal{C}^{\ast}\left(  Y\right)  $ the
complex of local chain currents and by $\mathcal{C}_{c\uparrow}^{\ast}\left(
Y\right)  $ the subcomplex with c/f support. These complexes compute cohomology with integer coefficients. Observe that since the stable manifolds are horned stratified, one can embed $_{\mathbb{Z}}\mathcal{S}_{c\uparrow}^{\ast} \subset \mathcal{C}_{c\uparrow}^{\ast} (Y)  $.\\

A chain current $R$ is transversal to an unstable manifold $U$ if the maps which locally define $R$ are transversal to $U$. The intersection number, which is an integer, is denoted by $ U \bullet R  =\lp U \wedge R \rp (\mathbf{1})$ (by definition, it vanishes if $dim R +  dim U \neq n $). We can now state the analogue of Lemma 2.3 for chain currents: the proof is omitted.
 
\begin{LE}\label{LE:extension-of-P}
The operators $\mathbf{T}$, $\mathbf{I}$ and $\mathbf{P}$ act on a local chain current $R$ on $Y$ with compact/forward support provided $R$ is transversal to each unstable manifold. The result is a local chain current with forward (resp. c/f) support. The Morse Chain Homotopy (for chain currents):
\begin{equation}\label{eq:MCHforS}
d ( \mathbf{T} R )+ \mathbf{T} (d R)= R -  \mathbf{P} (R)
\end{equation} 
holds as well as the expression:
\begin{equation}\label{eq:operatorPforS}
 \mathbf{P}(R)=\lim_{t \rightarrow +\infty}\phi_{t}^{\ast}(R)=\sum_{p\in Cr} \lp U_{p}\bullet R \rp  S_{p}
\end{equation}
In particular, $\mathbf{P}$ acts as the identity on any stable manifold $S_{p}$  and on the $\mathcal{S}$-complex as well. 
\end{LE} 

The last statement follows from either of the equalities in (\ref{eq:operatorPforS}). First because the stable manifolds are invariant under the flow and second because of the choice of orientations on the unstable manifolds.

\begin{RMK} Denote by $T(R)$ the set swept out by moving $R$ backward in time. Then this is the support of $\mathbf{T} (R)$ and one can picture the last equality in (\ref{eq:operatorPforS}) by visualizing dragging $R$ until it sticks on the critical points of the right index and clings to the corresponding stable manifold.
\end{RMK}

\begin{THE}\label{THE:Morse-integers} Suppose the flow is Weakly Proper and Smale. Then the maps 
\[
\mathbf{P}:\mathcal{C}_{c\uparrow}^{\ast}\left(  Y\right) \dashrightarrow \ _{\mathbb{Z}}\mathcal{S}_{c\uparrow}^{\ast}  \mbox{     and the inclusion    }\mathbf{I}:\  _{\mathbb{Z}}\mathcal{S}_{c\uparrow}^{\ast}
\hookrightarrow\mathcal{C}_{c\uparrow}^{\ast}\left(  Y\right)
\]
induce isomorphisms in cohomology
\[
H_{c\uparrow}^{k}\left(  Y,\mathbb{Z}\right) \approx H^{k} \left( \emph{{}}_{\mathbb{Z}}\mathcal{S}_{c\uparrow}^{\ast} \right)
\]
\end{THE}

The dashed map in the previous statement is defined only on chains transversal to the unstable manifolds (i.e. in the domain of $\mathbf{P}$ and $\mathbf{T}$).\\

\noi \textbf{Proof.} If the map $\mathbf{P}:\mathcal{C}_{c\uparrow}^{\ast}\left(  Y\right) \dashrightarrow \ _{\mathbb{Z}}\mathcal{S}_{c\uparrow}^{\ast}$ was defined on all of $\mathcal{C}_{c\uparrow}^{\ast}\left(  Y\right)$, the theorem would be a trivial consequence of Lemma \ref{LE:extension-of-P}. It hence suffices to show that in each class $[K]\in H^{p}_{c \uparrow}( \mathcal{C}^{\ast})$ there is a representative $[K']$ in the domain of $\mathbf{P}$ and $\mathbf{T}$.

Since the unstable manifolds make up a AB regular stratification of $Y$ (by Theorem \ref{THE:M}), the maps transversal to all the unstable manifolds are known to form a dense set among the possible ones. The result then follows by considering the pushforward of a 1-parameter deformation of $K$ $\Box$\\

The invariants extracted from the Morse function do not depend on the manifold alone (as is the case if $Y$ is compact); however, the next lemma proves that they do not change by bounded perturbations of the function.

\begin{LE}[Stability]\label{LE:stability} If $f_{0}$ and $f_{1}$ are two weakly proper functions on $Y$ whose difference is bounded (say by the constant $c\geq0$), then $f_{0}$ and $f_{1}$ determine the same family of compact/forward sets.
\end{LE}

\noi \textbf{Proof}. Note that $f_{1}^{-1}\left(  ]-\infty,a)\right)  \subset f_{0}^{-1}\left( ]-\infty,a+c)\right)  $ for any $a\in\mathbb{R}$, so that if $A$ is c/f with respect
to $f_{0}$ then $A$ is also c/f with respect to $f_{1}$ $\Box$\\

Actually, both the notions forward set and slab compact (cf. Definition \ref{DEF:forward}) are the same for $f_{0}$ and $f_{1}$.

\begin{RMK}[Morse inequalities] Suppose $f$ has only finitely many critical points, so that the $\mathcal{S}$-complex is finitely generated. The standard Morse inequalities are then a classical algebraic corollary of Theorem \ref{THE:Morse-integers}.
\end{RMK}

\begin{RMK}[$\Lambda$-Morse inequalities]\label{RMK:Lambda-inequalities} Later we will use the following more general result. Replace the integers $\Zs$ by any p.i.d. $\Lambda$ and suppose $(C^* , d) $ is any finitely generated complex of free $\Lambda$-modules. Then there are standard algebraic inequalities involving ranks and torsion numbers of the cohomology modules $H^k(C^*)$ and the ranks of the modules $C^k$.
\end{RMK}

\section{Circle Valued Morse Theory}
In this section, we present Novikov theory for circle valued maps as a variation of the previous Morse theory, governed by the addition of the action of the \emph{Novikov ring} on subsets of $Y$. With hindsight, our approach is very close in spirit to Novikov's original formulation in [N].\\

Suppose a circle valued Morse function $g:X\longrightarrow\mathbb{R}%
/\mathbb{Z}$ is given on the compact manifold $X$, and consider a gradient
vector field for $-g$, whose flow $\psi$ is Smale. Let $\sigma:\mathbb{R}\longrightarrow\mathbb{R}/\mathbb{Z}$ be the
quotient map and let
\[
\begin{array}{ccc}
Y & \overset{f}{\longrightarrow} & \mathbb{R}\\
\downarrow \rho &  & \downarrow\sigma\\
X & \overset{g}{\longrightarrow} & \mathbb{R}/\mathbb{Z}
\end{array}
\]
be the pullback covering. The group of deck transformations is isomorphic to the integers; we can choose a deck transformation $t:Y\longrightarrow Y$ as a generator of this group, so that the diffeomorphism $t$ and the Morse function $f$ are related by the equivariance (cf. later Remark \ref{RMK:generator}):
\begin{equation}\label{eq:equivariance}
f\left(  ty\right)  =f\left(  y\right)  +1  \, \, \, \text{ for any }y\in Y
\end{equation}

Using the covering map $\rho$, the gradient vector field and the flow $\psi$ can be lifted to a
vector field and a flow $\phi$ on $Y$. The flow $\phi$ is Smale and it's the gradient of the Morse function $-f$. In this cyclic covering case, $f$ is proper (not just weakly proper). In particular, the words ``forward'' and ``compact/forward'' have the same meaning; nevertheless we use the notation c/f, since it will be needed in the following section, dealing with the general case.

Observe that for each critical point $x$ of $g$, any arbitrary lifting $x_0$ is a critical point of $f$; all other critical points above $x$ can be written as $x_i = t^i (x_0)$ for $i \in \Zs$ and moreover $f(x_i) = f(x_0 ) +i$.\\ 

Consider now the \textbf{group ring }$\mathbb{Z} [
t,t^{-1} ]$ or $\mathbb{R} [  t,t^{-1} ]$ of the covering, i.e. the Laurent polynomials in $t$. Define the \textbf{Novikov ring} $\Lambda_{\mathbb{Z}}=\mathbb{Z}\left[
\left[  t\right]  \right]  [ t^{-1} ]$ and \textbf{Novikov field} $\Lambda_{\mathbb{R}}=\mathbb{R} \left[  \left[  t \right]  \right] [ t^{-1} ]$ to be the rings of formal Laurent series with finite principal parts. In the algebra literature, the ring $\Lambda_{\Zs}$ is referred to as the Laurent ring.

\begin{RMK}\label{RMK:generator} The choice of the generator $t$ for the group of deck transformations (rather than $-t$) determines the definition of the Novikov ring (using finite negative powers rather than finite positive powers).
\end{RMK}

The following algebraic facts are well known.

\begin{PROP}\label{PROP:algebra}
The ring $\Lambda_{\mathbb{Z}}$ is a $\mathbb{Z}  [  t,t^{-1} ]  $ module whereas  $\Lambda_{\mathbb{R}}$ is a $\mathbb{R} [  t,t^{-1} ]  $ module. Moreover:\newline
\vskip-2mm
Fact 1. The ring $\Lambda_{\mathbb{R}}$ is a field.\newline
\vskip-2mm
Fact 2. The ring $\Lambda_{\mathbb{Z}}$ is a principal ideal domain (pid).\newline
\vskip-2mm
Fact 3. The module $\Lambda_{\mathbb{Z}}$ is flat over $\mathbb{Z} [ t,t^{-1}] $. 
\end{PROP}
Briefly, the geometric series can be used to prove 1), and presenting $\Lambda_{\mathbb{Z}}$ as the localization at $\{t^{-1} \! , \! .., t^{-k} \!,\! .. \}$ of the ring of formal series $\mathbb{Z} [ [t]] $ provides a proof of 2) and 3). A detailed proof of the several variables case is given in the appendix.\\

As a consequence of the interaction between the deck map $t$ and $f$ (formula (\ref{eq:equivariance})), compact/forward sets can be defined using the deck map alone. This simple fact is the key for the compatibility of our Morse theory in the first part and the Novikov construction: 

\begin{LE}\label{LE:algebric-forward} A closed set $A\subset Y$ is a compact/forward set for $f$ if and only if there exists a compact set  $K\subset Y$ and an integer $N\in\mathbb{Z}$ such that $A\subset \bigcup_{n\geq N}t^{n}\left(  K\right)  $
\end{LE}

Let's now reconsider the complexes of forms and currents $\mathcal{E}%
_{c\uparrow}^{\ast}\left(  Y\right)  $, $\mathcal{D}_{c\uparrow}^{\prime\ast
}\left(  Y\right)  $, $\mathcal{C}_{c\uparrow}^{\ast}\left(  Y\right)  ,$ and
$\mathcal{S}_{c\uparrow}^{\ast}\left(  f\right)  $ discussed in the previous
section. 

First note that the deck map $t$ induces an action on currents (by pushforward) and that this action commutes up to sing with the differential $d$ on currents. Since the deck map $t$ also commutes with the flow $\phi$, Lemma \ref{LE:algebric-forward} implies that the induced action of $t$ is a self map of all the previous
complexes. 

The action of $t$ extends to an action of the group ring on each of the four complexes. Moreover, because of the support condition, this action extends to an action of the Novikov ring on those complexes.

\begin{PROP}\label{PROP:current Novikov complex}
The complex of currents ${}_{\mathbb{Z}}\mathcal{S}_{c\uparrow}^{\ast}\left(  f\right)$ is a free module over the Novikov ring $\Lambda_{\mathbb{Z}}$. The number of free generators in degree $k$ coincides with the number $N_k$ of critical points of index $k$ for the circle valued function $g$ (and hence is finite). 

Similarly, ${}_{\mathbb{R}}\mathcal{S}_{c\uparrow}^{k}\left(  f\right)$ is a vector space over $\Lambda_{\mathbb{R}}$ of dimension $N_k$.
\end{PROP}

\noi \textbf{Proof}. Any choice of a set of liftings $x\in Cr\left(  g\right)  \longmapsto x_0 \in Cr(f)$ for the critical points downstairs will provide a free basis for the $\Lambda_{\mathbb{Z}}$-module $\emph{{}}_{\mathbb{Z}}\mathcal{S}_{c\uparrow}^{\ast}\left( f\right)  $, consisting of the stable manifolds $S_{x_0} \in {}_{\mathbb{Z}}\mathcal{S}_{c\uparrow}^{\ast}$. Similarly for ${}_{\mathbb{R}}\mathcal{S}_{c\uparrow}^{\ast} $  $\Box$ \\
  
As a complex of $\Lambda_{\Zs}$-modules, $\lp {}_{\mathbb{Z}}\mathcal{S}_{c\uparrow}^{\ast}\left(  f\right) , d \rp$ is isomorphic to the so called \emph{Novikov complex}. As noted in the introduction, the latter is just a formal complex of $\Lambda_{\Zs}$-modules (see e.g. [P], [Sc2]). Our $\mathcal{S}$-complex realizes the Novikov complex geometrically as a subcomplex of currents, just as the current Morse complex geometrically realizes the Morse complex. By analogy, we will refer to $\lp {}_{\mathbb{Z}}\mathcal{S}_{c\uparrow}^{\ast}\left(  f\right) , d \rp$ as the \textbf{Current Novikov Complex}.\\

Recall the operators $\mathbf{T} :\mathcal{E}_{c\uparrow}^{\ast}\left(  Y\right)  \rightarrow
\mathcal{D}_{c\uparrow}^{\prime\ast}\left(  Y\right)$ and $\mathbf{P}: \mathcal{E}_{c\uparrow}^{\ast}\left(  Y\right)  \rightarrow {}_{\mathbb{R}}\mathcal{S}_{c\uparrow}^{\ast}\left(  f\right)
\subset \mathcal{D}_{c\uparrow}^{\prime\ast}\left(  Y\right)$. It is clear that they commute with the action of the diffeomorphism $t$, and so they are $\Lambda_{\mathbb{R}}$-linear maps.
Analogously, the inclusion map $\emph{{}}_{\mathbb{Z}}\mathcal{S}_{c\uparrow}^{\ast}\left(  f\right)
\hookrightarrow\mathcal{C}_{c\mathbb{\uparrow}}^{\ast}\left(  Y\right)  $ of stable manifolds into chain currents is a linear map of $\Lambda_{\mathbb{Z}}$-modules.

The realization of the formal Novikov complex as the current Novikov complex  $\mathcal{S}_{c\uparrow}^{\ast}$ produces the main result of this section, i.e. the characterization of Novikov homology as compact forward supported cohomology:    
\begin{THE}\label{THE:main} 
The map of $\Lambda_{\mathbb{R}}$-complexes $\displaystyle{\mathbf{P}:\mathcal{E}_{c\uparrow}^{\ast}\left(  Y\right)  \longrightarrow\emph{{}%
}_{\mathbb{R}}\mathcal{S}_{c\uparrow}^{\ast}\left(  f\right) }$ induces an isomorphism of finite dimensional $\Lambda_{\mathbb{R}}$-vector spaces
\[
H_{c\uparrow}^{k}\left(  Y,\mathbb{R}\right) = H^{k}\left(  \mathcal{E}^{\ast}_{c \uparrow} (Y) \right)  \approx H^{k}\left(  {}_{\mathbb{R}}\mathcal{S}_{c\uparrow}^{\ast}\left(  f\right)  \right) 
\]
Analogously, the inclusion map of $\Lambda_{\mathbb{Z}}$-complexes 
$\displaystyle{ {}_{\mathbb{Z}}\mathcal{S}_{c\uparrow}^{\ast}\left(  f\right)
\hookrightarrow\mathcal{C}_{\mathbb{\uparrow}}^{\ast}\left(  Y\right) }$
induces an isomorphism of finitely generated $\Lambda_{\mathbb{Z}}$-modules
\[
H^{k}\left(  {}_{\mathbb{Z}}\mathcal{S}_{c\uparrow}^{\ast}\left(
f\right)  \right)  \approx H^{k}\left( \mathcal{C}_{c\uparrow}^{\ast} \right) = H_{c\uparrow}^{k}\left(  Y, \mathbb{Z} \right) 
\]
\end{THE}

\noindent\textbf{Proof}. The isomorphism are a direct consequence of Theorems \ref{THE:Morse-reals} and \ref{THE:Morse-integers}, whereas the algebraic facts 1,2 imply that the cohomology spaces are finitely generated $\Box$ 

\begin{RMK}\label{RMK:Novikov inequalities} The $\Lambda_{\Zs}$-modules $H^{k}\left(  {}_{\mathbb{Z}}\mathcal{S}_{c\uparrow}^{\ast}\left( f\right)  \right) $ are thus isomorphic to Novikov homology, and since the Novikov ring $\Lambda_{\Zs}$ is a \emph{p.i.d.}, they have well defined ranks and torsion numbers, which are called \textbf{Novikov numbers} of the circle valued function $g$. 

The ``\emph{Novikov inequalities}'' relate those numbers to the number of critical points of the function $g$ (i.e. to the rank of the $\mathcal{S}$-complex) and are a purely algebraic corollary of the previous theorem, in analogy with the situation for Morse theory (see Remark \ref{RMK:Lambda-inequalities} or the discussion in [R], p.12).
\end{RMK}

Next, using our finite volume approach, we derive another characterization of Novikov homology, due to A. Pajitnov, [P3].\\

The sheaf cohomology groups $H_{cpt}^{k}\left( Y,\mathbb{Z}\right)  $ are standard homotopy invariants of $Y$ (isomorphic
to $H_{n-k}\left(  Y,\mathbb{Z}\right)  $, i.e. homology). The cohomology with c/f supports $H_{c\uparrow}^{\ast}\left(  Y,\mathbb{Z}\right)
$ depends not only on $Y$ but also on the family of compact/forward sets, though not on the particular function defining the c/f family, (nor on the Riemannian metric, a fortiori). We will modify the function in order to preserve the c/f supports and get a better gradient flow (and a better $\mathcal{S}$-complex), using the following lemma.

\begin{LE}\label{LE:Latour} Let $\omega$ be a closed one form on the compact manifold $X$. Then there exist a Riemannian metric on $X$, a Morse function $g_0 : X \rightarrow \Rs$, and a small constant $b>0$, such that the vector field  $\tilde{V} = V_0 +bV$, associated to the form $dg_0 + b \omega$, has nondegenerate singularities, fulfills the Smale condition, and its flow lines (even broken) are of uniformly bounded length. 
\end{LE}      

The lemma is probably valid (as the proof suggests) for a generic Riemannian metric and any Morse function $g_0$, provide $b$ is sufficiently small, but a self-contained proof would be long. We give a proof of the weaker case in the statement. See [L] for similar results.\\

\noi \textbf{Proof}. Start with any Morse function $g_0 : X \rightarrow \Rs$ with critical points distinct from those of $\omega$, and any Riemannian metric. If $b >0$ is sufficiently small, the critical points of $\tilde{V} = V_0 +bV$ are close to those of $V_0 = \nabla g_0$ and still non degenerate. Denote by $D = D(\varepsilon)$ the union of the $\varepsilon$-disks $D_{\varepsilon}(p)$ around the critical points $p$ of $\tilde{V}$.\\

\noi \textbf{Claim} \emph{For suitable $\varepsilon $ and $b$ the following three conditions hold}.
\vskip1ex
\noi $A$) \emph{The estimate $\| b V  \|  \leq \|  V_0 \|$ holds on $X -D $}.
\vskip1ex
\noi $B$) \emph{If $\tilde{\gamma}$ is a flow line of $\tilde{V}$ and $\tilde{\gamma}$ leaves one of the disks $D_{\varepsilon} (p)$, then it never comes back to it}.   
\vskip1ex 
\noi $C$) \emph{The flow of $\tilde{V}$ is Smale}.\\

First, we complete the proof assuming the claim. Note that since $\tilde{V} = V_0 +bV$, condition $A$) is equivalent to: 
\vskip1ex
\noi $A'$) \emph{The estimate $\| \tilde{V}  \|^2 \leq  2   V_0 \cdot \tilde{V}$ holds on $X -D $}.
\vskip1ex
This means that outside $D$ the values of $g_0$ grow along the flow lines of $\tilde{V}$. Since $X$ is compact and since $D$ contains the critical points of $\tilde{V}$, we can find a $C = C(D) >0$ so that 
\begin{equation}\label{eq:C}
\| \tilde{V}  \| \leq  C \|  \tilde{V} \|^2  \ \ \ \ \  \mbox{on } X-D
\end{equation}
Suppose $\tilde{\gamma}$ is any flow line of $\tilde{V}$. Then, by (\ref{eq:C}) and $A'$):
\begin{equation}\label{eq: abra}
 l(\tilde{\gamma} - D)    \leq \,     C \int_{\tilde{\gamma} - D} \| \tilde{V}  \|^2 dt 
 \leq  2 C \int_{\tilde{\gamma} - D}  V_0 \cdot  \tilde{V}   dt = 2 C \int_{\tilde{\gamma} -D} dg_0 
\end{equation}
Because of this estimate, the length of any single arc of $\tilde{\gamma} - D$ is uniformly bounded. By condition $B$), the number of such arcs cannot exceed the number of critical points. The number of arcs in $\tilde{\gamma} \cap D$ is then also less then the number of critical points, and the length of each such arc is obviously uniformly bounded. The lemma follows.

It remains to prove the claim. Choose a small local coordinate patch $U$ near each critical point of $g_0$, so that the form $\omega$ is the differential of a coordinate in each such $U$. By cutting and smoothing $g_0$ near its critical points (and shrinking the $U$'s), we can assume that $g_0$ is a quadratic form in standard form, in the chosen coordinates. Finally, assuming the metric is the standard metric in the coordinates of $U$, we can write:
\beq \label{eq:Lx}
V_0 (x) = Lx  \hspace{1ex} \mbox{ and } \hspace{1ex} V (x) = v \hspace{4ex} \mbox { on } U
\eeq
where $v$ is a constant vector and $L$ a diagonal matrix with eigenvalues $\pm1$. 
 
If $b >0$ is sufficiently small, the critical points of $\tilde{V} = V_0 +bV$ are $w = -L^{- \! 1}(bv)$, one in each chart $U$ as above. Setting $y = x - w$:
\begin{equation}\label{eq:Ly}
\tilde{V} (y) = L y \hspace{4ex} \mbox { on } U
\end{equation}
and hence the critical points of $\tilde{V}$ are nondegenerate.

Let $D_{\varepsilon} (w)$ be the product $\varepsilon$-disk in the (split) $y$-coordinates above, centered at the critical point $w = -L^{- \! 1}( bv) \in U$. One can choose $\varepsilon$ small enough so that $D_{(k+1)\varepsilon} (w) \subset U$ and $k\! > \! 2$. In this case, since the flow is linear on $U$, if a flow line leaves $D_{\varepsilon} (w)$, then it necessarily also leaves $D_{(k+1)\varepsilon} (w)$. Moreover, setting $\tilde{\omega} = dg_0 + b\omega$ (the 1-form corresponding to $\tilde{V}$), one can estimate, over any connected arc of $\tilde{\gamma} \cap ( D_{(k+1)\varepsilon} -  D_{\varepsilon})$ leaving $D_{\varepsilon}$:
\beq \label{eq:variation}
\int\limits_{\tilde{\gamma} \cap ( D_{(k+1)\varepsilon} -  D_{\varepsilon} ) } \tilde{\omega} \approx ((k+1)^2 \varepsilon ^2 - \varepsilon^2 ) - (\varepsilon^2) > 4 \varepsilon^2 \approx 2 ( \max\limits_{D_{\varepsilon}} (\tilde{g})- \min\limits_{D_{\varepsilon}}(\tilde{g}) ) 
\eeq
Fixing $\varepsilon$, we can choose $b$ sufficiently small (then $\tilde{\omega} \approx dg_0$) so that condition $A$) is true and the following holds, because of (\ref{eq:variation}): 
\beq\label{eq:cadabra} 
\int\limits_{\tilde{\gamma} \cap \lp D_{(k+1)\varepsilon} -  D_{\varepsilon}\rp } dg_0  > \max\limits_{D_{\varepsilon}} (g_0)- \min\limits_{D_{\varepsilon}}( g_0) 
\eeq
Since a flow line leaving $D_{\varepsilon} (w)$, it necessarily also leaves $D_{(k+1)\varepsilon} (w)$, the estimate (\ref{eq:cadabra}) implies that such a flow line cannot return to $D_{\varepsilon} (w)$ again. This proves conditions $A$) and $B$) in the claim. Condition $C$) can finally be obtained by slightly perturbing the Riemannian metric (and hence $\tilde{V}$) on the sets $U- D_{(k+1)\varepsilon}$, with no influence on the behavior established above $\Box$ \\

Back to Novikov theory, we apply the previous lemma to the case when $\omega = dg$ (our circle valued function), obtaining a Riemannian metric on $X$, a Morse function $g_0: X \rightarrow \Rs$ and a constant $b \ll 1$. The lift of $g_0$ to $Y$ is a bounded Morse function $f_0 : Y \rightarrow \Rs$. The function $\tilde{f} = f_0 +b f $, i.e. the lifting of $\tilde{g} = g_0 +bg$, is then proper, being a bounded perturbation of the proper function $f$. Moreover, the previous lemma implies that $\tilde{f}$ is Morse and that there is a uniform bound on the length of the flow lines of the gradient of $\tilde{f}$ (for the lifted metric). By rescaling, the same is true for the flow lines of the gradient of $f + c  f_0$, where $1 /b = c \gg 0$.\\  

Using the new Morse function $f + c  f_0 $ and the new Riemannian metric on $Y$, one can construct a new $\mathcal{S}$-complex. The Stability Lemma \ref{LE:stability} proves that the two functions, $f$ and $f + c  f_0 $, define the same family of compact/forward sets, since $f_0$ is bounded. The cohomologies of the two $\mathcal{S}$-complexes are thus isomorphic, being an invariant of $Y$ and the family of c/f sets only. 

In the rest of the section, we will hence replace the old function $f$ (and Riemannian metric) with the new function $f + c  f_0 $ (and the new metric): this is sometimes called \emph{Latour's trick}. The advantage is that we can now assume a uniform bound on the length of any (even broken) flow line for $f$.\\
        
The isomorphism of $\Lambda$-modules $H^{k}\left(  \emph{{}}_{\mathbb{Z}%
}\mathcal{S}_{c\uparrow}^{\ast}\left(  f\right)  \right)  \approx
H_{c\uparrow}^{k}\left(  Y,\mathbb{Z}\right)  $ remains valid. Also, each stable manifold $S_{x}$ is now relatively compact in $Y$, therefore $\left[  S_{x}\right]  $ has compact support and its boundary consists of a finite sum of other stable manifolds. In particular the space
$\emph{{}}_{\mathbb{Z}}\mathcal{S}_{cpt}^{\ast}\left(  f\right)  $, made up of finite sums of stable manifolds, is closed under taking boundary, i.e. it is a complex. Moreover, the operator $\mathbf{P}$ now maps $\mathcal{E}_{cpt}^{\ast}\left(
Y\right)  $ to $\emph{{}}_{\mathbb{R}}\mathcal{S}_{cpt}^{\ast}\left(
f\right)  $, and  $\mathbf{T}$ is a chain homotopy between $\mathbf{P}$ and the inclusion $\mathbf{I}:\mathcal{E}_{cpt}^{\ast}\left(  Y\right)  \hookrightarrow
\mathcal{D}_{cpt}^{\prime \ast}\left(  Y\right)  $.

Consequently, there are isomorphisms of real vector spaces and abelian groups:%
\[
H_{cpt}^{k}\left(  Y,\mathbb{R}\right)  \approx H^{k}\left(
\mathcal{E}_{cpt}^{\ast}\left(  Y\right)  \right)  \approx H^{k}\left(
\emph{{}}_{\mathbb{R}}\mathcal{S}_{cpt}^{\ast}\left(  f\right)  \right)
\]%
\[
H_{cpt}^{k}\left(  Y,\mathbb{Z}\right)  \approx H^{k}\left(
\mathcal{C}_{cpt}^{\ast}\left(  Y\right)  \right)  \approx H^{k}\left(
\emph{{}}_{\mathbb{Z}}\mathcal{S}_{cpt}^{\ast}\left(  f\right)  \right)
\]

The (covering) group ring $\mathbb{Z}\left[  \pi\right]  =\mathbb{Z}\left[
t,t^{-1}\right]  $ of Laurent polynomials acts on $\emph{{}}_{\mathbb{Z}%
}\mathcal{S}_{cpt}^{\ast}\left(  f\right)  $ and $\mathcal{C}_{cpt}^{\ast
}\left(  Y\right)  $. Therefore $\emph{{}}_{\mathbb{Z}}\mathcal{S}_{cpt}%
^{\ast}\left(  f\right)  \underset{\mathbb{Z}\left[  \pi\right]  }{\otimes
}\Lambda_{\mathbb{Z}}$ and $\mathcal{C}_{cpt}^{\ast}\left(  Y\right)
\underset{\mathbb{Z}\left[  \pi\right]  }{\otimes}\Lambda_{\mathbb{Z}}$ are
complexes of $\Lambda_{\mathbb{Z}}$ modules and there are isomorphisms of
$\Lambda_{\mathbb{Z}}$ modules:%
\[
_{\mathbb{Z}}\mathcal{S}_{cpt}^{\ast}\left(  f\right)  \underset
{\mathbb{Z}\left[  \pi\right]  }{\otimes}\Lambda_{\mathbb{Z}}={}_{\mathbb{Z}%
}\mathcal{S}_{c\uparrow}^{\ast}\left(  f\right)  \text{ \ \ and \ \ }%
\mathcal{C}_{cpt}^{\ast}\left(  Y\right)  \underset{\mathbb{Z}\left[
\pi\right]  }{\otimes}\Lambda_{\mathbb{Z}}=\mathcal{C}_{c\uparrow}^{\ast
}\left(  Y\right)
\]

Taking cohomology of the complexes, and using the algebraic Fact 3, the previous discussion proves the following characterization, cf.[P3]:

\begin{THE}[Pajitnov]\label{THE:pajitnov} As finitely generated $\Lambda_{\mathbb{Z}%
}$-modules, Novikov homology is isomorphic to compactly supported cohomology with coefficients in the Novikov ring $\Lambda_{\Zs}$. That is:
\[
H^{k}\left( {}_{\Zs}\mathcal{S}^*_{c \uparrow} \right) \approx H_{cpt}^{k}\left(  Y , \mathbb{Z}\right)  \underset{\mathbb{Z}\left[
\pi\right]  }{\otimes}\Lambda_{\mathbb{Z}}
\]
\end{THE}

Pajitnov's proof of this theorem is just for the case of cyclic coverings, discussed in the present section. The more general case of irrational 1-forms, discussed in the next section, can be obtained by approximating irrational forms by rational forms: a rather detailed proof is due to Schutz [Sc], §4.2. Our proof in the general case does not differ from the proof for cyclic coverings presented above.   

\section{Novikov Theory: General Case}\label{sec:Novikov-general}

In this section we consider the general case of Novikov theory, which deals with closed Morse 1-forms on compact manifolds. Contrary to the previous section, the Morse function arising in this case is not proper. Our finite volume approach to Morse theory is thus particularly fruitful for the geometric interpretation of the general Novikov theory.\\

Let $\omega$ be a Morse $1$-form on the compact Riemannian manifold $X$,
i.e. a closed one-form with nondegenerate singularities. Its gradient vector field determines a flow $\psi$ on $X$. Using $\psi$, one can define global stable and unstable manifolds. We will assume the flow to be Smale (i.e.
all stable and unstable manifolds have to intersect transversally: this is
known to be a generic condition for this kind of gradient field).

Let $q-1$ be the irrationality index of $\omega$ and $\chi=\left( \chi_{1},\ldots ,\chi_{q} \right) \in \Rs^q $ denote its periods: this means that one can write $\omega = \sum \chi_i \omega^i $ where the closed forms $\omega^i$ define classes in $H^1(X,\mathbb{Z})$ and $q$ is the minimal number of such forms. One can assume the periods to be positive numbers. 

Let $\rho:Y\rightarrow X$ be a minimal covering such that
$\omega$ pulls back to an exact form, say $df$, with $f:Y\longrightarrow
\mathbb{R}$. The group $\pi$ of deck transformations of $\left(  Y,\rho\right)  $
is a free abelian group with $q$ generators, say $t_{1},..,t_{q}$ (i.e.
$\pi\approx\mathbb{Z}^{q}$) and the \textbf{group ring } is the ring of Laurent polynomials in $q$ variables:
\[
\mathbb{Z}\left[  \pi\right] =\mathbb{Z}\left[  t, t^{-1} \right] \text{ \ \ \ or \ \ \ } \mathbb{R}\left[
\pi\right]  =\mathbb{R}\left[  t,t^{-1}\right] \text{ \ \ \ \  ,\  \ \ \ } t= (t_{1},..,t_{q}) 
\]

Fixing the choice of generators, we can assume that the equivariance relations $f\left(  t_{i}\left(  y\right)  \right)  =f\left( y\right)  +\chi_{i}$ hold for each $i=1,..,q$. As noted in Remark \ref{RMK:generator}, other choices would make necessary a different definition of the Novikov Ring (see Definition \ref{DEF:Novikov-Ring}). 

Note that if $q=1$ the covering is cyclic and the one form $\omega$ can be seen as the differential of a circle valued function, which is exactly the case in the previous section. \\

To define the appropriate Novikov ring in several variables, let the vector of periods
$\chi$ also denote the linear form on $\mathbb{R}^{q}$ defined by
$\chi\left(  v\right)  =\chi\mathbf{\cdot}v$. A subset of the lattice $F\subset\mathbb{Z}^{q}$ is called:

\bigskip

1) \textbf{slab compact} \emph{if }$F$\emph{ intersected with each slab }$\chi
^{-1}\left(  \left[  a,b\right]  \right)  $\emph{ is compact (i.e. finite),}

2) \textbf{forward}\emph{ if }$F\subset\chi^{-1}\left(  [a,+\infty)\right)
$\emph{ for some }$a\in\mathbb{R}$,

3) \textbf{compact/forward}\emph{ or c/f if }$F$\emph{ is both slab compact and forward.}

\bigskip

Consider now the formal Laurent series in $q$ variables $\alpha=\sum a_{n}t^{n}$ , where
$t=\left(  t_{1},..,t_{q}\right)  $ and $n=\left(  n_{1},..,n_{q}\right)  $.
The \textbf{support} of $\alpha$, denoted by spt$ \lp \alpha \rp  $,
 consists of all $n\in\mathbb{Z}^{q}$ such that $a_{n}\neq 0$.

\begin{DEF}\label{DEF:Novikov-Ring} The \textbf{Novikov ring} $\Lambda_{\Zs}$ (defined by the linear form $\chi : \mathbb{R}^{q} \rightarrow \Rs$) consists of all formal Laurent series $\alpha=\sum a_{n}t^{n}$ with compact/forward support (with respect to $\chi$) and with integer coefficients $a_{n}\in\mathbb{Z}$.
Similarly, one can choose real coefficients to define $\Lambda_{\mathbb{R}}$.
\end{DEF}

The ring structure is defined as follows. The support of $\alpha$ is co\-mp\-act/for\-ward if and only if spt$\lp  \alpha \rp \cap\chi^{-1}\left(
(-\infty,a]\right)  $ is finite for all $a\in\mathbb{R}$. Consequently, given $\alpha,\beta\in\Lambda$, the
Cauchy product $\gamma=\alpha\beta$ is defined by the finite sums $c_{n}%
=\sum\limits_{p+q=n}a_{p}b_{q}$, and $\gamma$ has c/f support too. The following proposition, due to Sikorav and Pajitnov, is proved in the appendix.

\begin{PROP}\label{PROP:algebra2} The algebraic Facts 1,2,3 of Prop. \ref{PROP:algebra} hold true for the Novikov ring.
\end{PROP}

Let's now return to the covering $\rho:Y\rightarrow X$ which ``unties'' the Novikov form $\omega$. As it has been done in the case of cyclic coverings, the
gradient vector field and the flow $\psi$ can be lifted to a vector field and
flow $\phi$ on $Y.$ The flow $\phi$ is the gradient of the Morse
function $f$ and is again Smale, though $f$ is \textbf{not} proper now. Of course, upstairs (i.e. on $Y$) there are no closed orbits and the critical points of $f$ are just the preimages of the critical points of $\omega$.

\begin{LE} The lifted flow $\phi$ is weakly proper.
\end{LE}

The lemma and the proof below hold for any covering where $\omega$ pulls back to an exact form. In particular the lifted flow is weakly proper for the universal covering.\\

\noindent\textbf{Proof}. Suppose $\bar{\gamma}:[0,+\infty\lbrack\rightarrow Y$
is a forward flow-half line of $\phi$ which is not relatively compact in $Y$
(i.e. $\bar{\gamma}$ does not converge to a critical point); we just need to
show that $f$ is unbounded on $\bar{\gamma}$. Consider the projected curve
$\gamma=\rho\left(  \bar{\gamma}\right)  $, which is a forward flow-half line
for $\psi$. Since $\rho^{\ast}\left(  \omega\right)  =df$ , and
\[
\lim\limits_{s\rightarrow+\infty}\left[  f\left(  \bar{\gamma}\left(
s\right)  \right)  -f\left(  \bar{\gamma}\left(  0\right)  \right)  \right]
=\lim\limits_{s\rightarrow+\infty}\int_{0}^{s}\frac{d}{dt}f\left(  \bar
{\gamma}\left(  t\right)  \right)  dt=\int_{\bar{\gamma}}df=\int_{\gamma
}\omega
\]
we are left to prove that $\int_{\gamma}\omega=+\infty$. Observe that $\gamma$
cannot converge to a critical point for $\psi$ in $X$, otherwise $\bar{\gamma
}$ would also converge to a critical point for $\phi$ in $Y$. Moreover, for
any open set $D\subset X$ containing all the critical points, there exists a
constant $c>0$, determined by $\left|  \omega\right|  $ on $X\backslash D$,
such that for any piece of an integral curve $\alpha$ contained in
$X\backslash D$, the estimate $\int_{\alpha}\omega=\int\left|  \dot{\alpha
}\right|  ^{2}dt>c\int\left|  \dot{\alpha}\right|  dt$ holds. Since $\gamma$
doesn't converge to a critical point, we can choose $D$ so that $\gamma$ has
unbounded length in $X\backslash D$ $\Box$\\

The lemma above permits one to apply our non compact Morse theory to the Novikov setting. On the other hand, the next lemma allows one to define compact/forward sets using only the covering group $\pi$ and not the function $f$. These two results are the basis for our development.       

\begin{LE} A closed set $A\subset Y$ is a
compact/forward set if and only if there exists a compact set $K\subset
Y$ and a $c/f$ set $F$ in the lattice $\mathbb{Z}^{q}$ such that $A$ is contained in the union of the sets $t^{n}\left(  K\right)  $ over $n\in F$ .
\end{LE}

Because of this lemma, the complexes $\mathcal{S}^*_{c \uparrow}$, $\mathcal{E}^*_{c \uparrow}$ and $\mathcal{C}^*_{c \uparrow}$ are naturally $\Lambda$-modules, the action being induced by pushforward by the diffeomorphisms $t_i$. In particular, $\mathcal{S}^*_{c \uparrow}$ is free, and finitely generated over $\Lambda_{\Zs}$ since $X$ is compact. As in the cyclic covering case, the complex $\mathcal{S}^*_{c \uparrow}$ is a geometric realization of the Novikov complex. 

Now we can state the main theorem, extending word for word the results in the cyclic covering case:

\begin{THE}\label{THE:main-general} Theorems \ref{THE:main} and \ref{THE:pajitnov} remain true for the general case $q \geq 1$.
\end{THE}

The proofs are identical to the cyclic covering case. Only note that to apply Latour's trick in the proof of Theorem \ref{THE:pajitnov}, one needs to add a multiple of an exact form (instead of a function). 

Exactly as before, the Novikov inequalities follow as an algebraic corollary, cf. Remark \ref{RMK:Novikov inequalities}. We also point out:

\begin{RMK}[Topological Stability] Any two Morse 1-forms in the
same cohomology class in $H^{1}\left(  X,\mathbb{R}\right)  $ define the same
c/f sets on $Y$. In fact they differ by the differential of a bounded function,
since $X$ is compact. Their liftings to $Y$ thus differ by a bounded
function and the Stability Lemma \ref{LE:stability} applies.
\end{RMK}

We end this section by discussing a slight variation of the previous theory.\\

 \noi \textbf{ The Novikov conical ring}\\

In the paper [N2], Novikov used a different ring of formal series, instead of (what we call) the Novikov ring, to develop his theory. We introduce it in the appendix as the ``Novikov conical ring'' $ {}_{_<}\!\Lambda$; it is a subring of the Novikov ring $\Lambda$. 

One can define the Novikov complex in a completely similar manner using the Novikov conical ring $ {}_{_<}\!\Lambda$. This produces a second complex, made up of free modules over $ {}_{_<}\!\Lambda$, again generated by the critical points of the Morse 1-form. Although algebraically distinct from the previous Novikov complex (over $\Lambda$), the second one (over $ {}_{_<}\!\Lambda$) has the same homology properties, because the subring $ {}_{_<}\!\Lambda$ is a flat $\Lambda$-algebra, as proved in the appendix. This implies, for example, that the Novikov numbers are the same, so that the Novikov inequalities are also the same.\\

With our approach, it is possible to first compute the cohomology of the Novikov complex as sheaf cohomology with compact forward supports in $Y$, and later, at the end, use the module structure over the Novikov ring. The flexibility of the finite volume technique allows one to follow a similar path to describe the cohomology of the ``conical Novikov complex''.

In fact, Lemma \ref{LE:algebric-forward} can be used as a definition for the family of \textbf{conical forward} sets, replacing $\Lambda$ by $ {}_{_<}\!\Lambda$. One can then define the complexes of forms, currents and the $\mathcal{S}$-complex with supports in the family of conical forward sets, and show that the Morse Chain Homotopy in Section 2 still holds; this just involves checking supports. 

We already used a similar argument (i.e. define an $\mathcal{S}$-complex over a family of supports different from the compact/forward ones and study the MCH equation with such supports), during the proof of Pajitnov's Theorem in Section 3. 

The cohomology of the Novikov conical complex is the sheaf cohomology of $Y$ with supports in the family of conical forward set, and is endowed with a natural module structure over the Novikov conical ring. The interested reader might complete the proof, or see [M].

\section{Backward-Forward Dualities}

Throughout this section, we will be in the geometric framework of Section \ref{sec:Novikov-general} (Novikov theory, $q \geq 1$), adopting the same notations. The aim is to study duality for Novikov homology, both as $\Lambda$-modules and as real (topological) vector spaces.\\

Let's start by recalling that a closed set $B\subset Y$ is \emph{backward} iff $B \subset f^{-1}(]-\infty , b]) $ for some $b \in \Rs$. A \emph{compact/backward} set is a backward set which is also slab compact (cf. Definition \ref{DEF:forward}). 

The arguments in Sections 1,2 can be repeated for the operator
\[
\mathbf{R} = \lim_{t \rightarrow - \infty} \phi_t^*
\]
obtaining completely analogous results for cohomology with compact backward supports (and for backward cohomology as well, recall Remark \ref{RMK:forward}). 

One defines the two \emph{$\mathcal{U}$-complexes} $\mathcal{U}^{\ast}_{\downarrow}(f)$ and $\mathcal{U}^{\ast}_{c \downarrow}(f)$ as made up of currents of the form $\sum_{x\in B}a_x \left[ U_x \right] $, where $B$ is a backward (resp. c/b) set of critical points, the $a_x$ are numbers (reals or integers) and the $[U_x]$ is the current of integration over the unstable manifold of the critical point $x$. For the differential on these two complexes it's convenient to choose the current boundary $\partial$. 

The operator $\mathbf{R}$ satisfies a chain homotopy similar to the Morse chain homotopy MCH, and the Theorems \ref{THE:Morse-reals} and \ref{THE:Morse-integers} still hold for backward and compact backward supports, using the $\mathcal{U}$-complexes and the operator $\mathbf{R}$ in place of $\mathcal{S}$ and $\mathbf{P}$.\\  
    
Like the $\mathcal{S}$-complexes, $\mathcal{U}^{\ast}_{\downarrow}(f)$ and $\mathcal{U}^{\ast}_{c \downarrow}(f)$ also have natural structures of modules over the Novikov ring $\Lambda$. We let the deck transformations $t_i$ act by pullback on the unstable manifolds, and extend the action by linearity to $\Lambda$. This action is different from the one on the $\mathcal{S}$-complex, which is via pushforward.

It is important to observe that $\mathcal{U}^{\ast}_{c \downarrow}(f)$ is finitely generated over $\Lambda$, whereas $\mathcal{U}^{\ast}_{\downarrow}(f)$ is not (if $q\! > \! 1$, otherwise they coincide). In both cases, though, the differential $d$ is a $\Lambda$-linear map.\\

\textbf{Lambda Duality} 

\bigskip

In this subsection $\Lambda$ denotes the Novikov ring with integer coefficients. The $\mathcal{S}$ and $\mathcal{U}$-complexes are intended with integer coefficients as well.\\

 Any critical point $x\in X$ of index $k$ determines a rank one $\Lambda$-submodule $\mathcal{S}^k_x \subset \mathcal{S}_{c\uparrow}^k$, generated by $S_{x_0}$, where $x_0 \in Y$ is any fixed lifting of $x$. Clearly $\mathcal{S}^k_x$ does not depend on the choice of $x_0$; on the other hand any $S \in \mathcal{S}^k_x$ determines a unique $\lambda \in \Lambda$ (which indeed depends on $x_0$) for which $S = \lambda S_{x_0}$. 

Similarly, the critical point $x$ determines a rank-one $\Lambda$-submodule $\mathcal{U}_x^{n-k} \! \subset \mathcal{U}_{c \downarrow }^{n-k} $, generated by $U_{x_0}$; as before, $U\in\mathcal{U}_{x}^{n-k}$ uniquely determines $\mu\in\Lambda$ with $U=\mu U_{x_{0}}$.\\           

Given $U\in\mathcal{U}_{x}^{n-k}$ and $S\in\mathcal{S}_{y}^{k}$, define the $\Lambda
$-pairing $U \bullet_{_{\Lambda}} S$ to vanish if $x \neq y$ and to be
\[
U \bullet_{_{\Lambda}} S =\lambda \mu
\]
when $x= y$, $x_0$ is any fixed lifting of $x$, and $\lambda$ and $\mu$ are determined by the relations $U=\mu U_{x_{0}}$ and $S=\lambda S_{x_{0}}$. 

\begin{PROP-DEF}
The $\Lambda$-pairing $\bullet_{_{\Lambda}} $ is independent of the choice of the critical point $x_{0}\in Y$ lifting $x\in Cr\left(  \omega\right) \subset X$, and extends by linearity to a $\Lambda$-bilinear pairing $\mathcal{U}_{c\downarrow}^{\ast}(Y) \times \mathcal{S}_{c\uparrow}^{\ast} (Y) \rightarrow \Lambda$. 
\end{PROP-DEF}

\noindent \textbf{Proof}. Suppose $x_{0}^{\prime} = t^{n}x_{0}$ is another choice of critical point above $x$. Then $S_{x_{0}^{\prime}}=S_{t^{n}x_{0}}%
=t^{n}S_{x_{0}}$ while $U_{x_{0}^{\prime}}=U_{t^{n}x_{0}}=t^{-n}U_{x_{0}}$ for the same $n\in \Zs^q$ so that $\lambda=\lambda^{\prime}t^{n}$ and $\mu
=\mu^{\prime}t^{-n}$ and $\lambda\mu=\lambda^{\prime}\mu^{\prime}$ $\Box$\\

For any critical point $x \in X$, and an arbitrary lifting $x_0 \in Y$, we choose an orientation on the stable manifold $S_{x_0}$ and orient all other $S_{x_i}$ via the covering maps $t^i$ ($i$ is a multiindex here). This induced an orientation on the unstable manifolds $U_{x_i}$ (since $Y$ is oriented) so that the intersection number $U_{x_0} \bullet S_{x_0} = \lp U_{x_0} \wedge S_{x_0}\rp  \lp \mathbf{1} \rp = 1$ and the intersection of $U_{x_0}$ with any other stable manifold vanishes. 

\begin{RMK} The constant term in the series $U \bullet_{_{\Lambda}} \! S \in \Lambda$ coincides with the intersection number $ U \bullet S  \in \mathbb{Z}$. More generally, the coefficient of the monomial $t^i$ in $U \bullet_{_{\Lambda}} \! S $ is given by $(t^{-i} U) \bullet S = U \bullet (t^{-i}S) $. 
\end{RMK}

The pairing $ \bullet_{_{\Lambda}}$ allows one to identify $\mathcal{U}_{c \downarrow}^{\ast} \approx Hom_{\Lambda}  \lp \mathcal{S}_{c \uparrow}^{\ast} , \Lambda \rp $ as $\Lambda$-modules. Recall that $\partial$ denotes the current boundary (which agrees with the current differential $d$ up to a sign depending on the degree of the current).

\begin{LE} The following complexes of finitely generated $\Lambda$-modules are isomorphic:
\[ 
\lp \mathcal{U}_{c \downarrow}^{\ast}(Y), \partial  \rp  \approx \text{Hom}\lp  \lp \mathcal{S}_{c \uparrow}^{\ast}(Y) , d \rp \, , \Lambda \rp
\] 
\end{LE}

\noi \textbf{Proof}. The only thing to prove is that the operator $\partial$ on $ \mathcal{U}_{c \downarrow}^{\ast}$ is the adjoint via $ \bullet_{_{\Lambda}}$ of the operator $d$ on $\mathcal{S}_{c \uparrow}^{\ast} $. This is a trivial consequence of the previous remark and of the fact that $d$ and $\partial$ are adjoint for the current intersection pairing between the $\mathcal{U}$ and $\mathcal{S}$-complex (cf. later, Proposition \ref{PROP:duality-U-S-complexes}) $\Box$\\
  
The next theorem is a standard algebraic consequence of this lemma. A weaker result was obtained by Pajitnov in [P2] from a different point of view, cf. also [Sc2]. 

\begin{THE}[Lambda Duality]\label{THE: Lambda duality} The $\Lambda$-pairing between the $\mathcal{U}_{c \downarrow}$ and $\mathcal{S}_{c\uparrow }$ complexes induces a $\Lambda$ bilinear form $\wedge : H^k_{c \downarrow} (Y, \Zs) \times H^{n-k}_{c \uparrow} (Y, \Zs ) \rightarrow \Lambda $ with the following two properties.
\vskip2ex

\noi \hskip1ex $1$. A class $u \in H^k_{c \downarrow} (Y, \Zs)$ satisfies $u \wedge s =0$ for all $s\in H^{n-k}_{c \uparrow} (Y, \Zs )$ iff $u$ is a $\Lambda$-torsion class.
\bigskip

Let both $u \in H^k_{c \downarrow} (Y, \Zs)$ and $s\in H^{n-k+1}_{c \uparrow} (Y, \Zs )$ be $\Lambda$-torsion classes. Then the ``\emph{linking number}'' $ \ell (u,s)= \frac{1}{\lambda} V \bullet_{_{\Lambda}} S \in \Lambda_{\Rs}/\Lambda $ is well defined where  $U\in u$, $S \in s$ and $d V = \lambda U$ for some $\lambda \in \Lambda$. 
\vskip2ex

\noi \hskip1ex $2$. Suppose $u\in H^k_{c \downarrow} (Y, \Zs)$ is a $\Lambda$-torsion class. Then $\ell (u,s)$ vanishes for all $\Lambda$-torsion classes $s \in H^{n-k-1}_{c \uparrow} (Y, \Zs )$ iff $u=0$.
\end{THE} 

\vskip4mm

\textbf{Topological Vector Space Duality}

\bigskip

In this subsection, real coefficients are understood for the Novikov ring $\Lambda$ (which is a field) and for the $\mathcal{S}$ and $\mathcal{U}$-complexes.\\

We are going to consider the Novikov homology $H^k_{c \uparrow} (Y, \Rs)$ as a real topological vector space. Quite a lot of the following discussion might as well be done for the case of an arbitrary Morse function on an arbitrary manifold $Y$. However, the geometric framework of Novikov theory allows simplifications of the proofs, and we will just treat this case here. \\

We will see that the complexes of forms, currents and the $\mathcal{S}$ and $\mathcal{U}$ complexes (with suitable supports) bear natural topologies. It follows that the cohomology of each complex carries a topology too. Contrary to the Lambda duality, in order to get a duality over the reals we need to pair ``compact/forward'' objects with ``backward'' ones. Our aim is to prove the following.  

\begin{THE}\label{THE:top.isomorphism} The algebraic isomorphism stated in Theorem \ref{THE:Morse-reals}
\begin{equation}\label{eq:THE:top.isomorphism}
H^{k}\left( \mathcal{E}_{c \uparrow}^{\ast} (Y)\right) \overset{\mathbf{P}}{\approx} H^{k}\left( \mathcal{S}_{ c\uparrow}^{\ast} (Y)\right) \overset{\mathbf{I}}{\approx}  H^{k}\left( \mathcal{D'}_{c \uparrow}^{\ast} (Y)\right) 
\end{equation}
are topological isomorphism among Hausdorff locally convex spaces. The same is true for the analogous isomorphisms:
 
\begin{equation}\label{eq:THE:top.isomorphism2}
H^{k}\left( \mathcal{E}_{\downarrow}^{\ast}(Y) \right) \overset{\mathbf{R}}{\approx} H^{k}\left( \mathcal{U}_{\downarrow}^{\ast} (Y)\right) \overset{\mathbf{I}}{\approx}  H^{k}\left( \mathcal{D'}_{\downarrow}^{\ast}(Y) \right)
\end{equation}

Moreover, the first spaces are in duality with the second ones and the pairing is equivalently  induced by the pairing between c/f forms and backward currents or between c/f currents and backward forms or by the current intersection between  $\mathcal{S}_{c \uparrow}^{\ast}$ and $\mathcal{U}_{\downarrow}^{\ast}$.
\end{THE}

We postpone the proof until after some preliminaries. Since all the complexes are considered on the covering manifold $Y$ (the base $X$ is not used below), we will sometimes drop the explicit reference to it.\\

First, we fix the topologies on the spaces involved, starting with the spaces of smooth forms. 

For any closed set $A\subset Y$ let $\mathcal{E}_{A}^{\ast} \subset\mathcal{E}^{\ast} (Y) $ denote the space of smooth forms with support in $A$. This is a closed subspace of the Frechet space $\mathcal{E}^{\ast} (Y) $ and is hence Frechet. Since algebraically $\mathcal{E}_{c \uparrow}^{\ast} (Y) = \lim_{\rightarrow} \mathcal{E}_{A}^{\ast}$, where the direct limit is taken over the directed family of c/f sets, it's natural to give the locally convex inductive limit topology to $\mathcal{E}_{c \uparrow}^{\ast} (Y)$. Similarly, we can give the inductive limit topology to $\mathcal{E}_{\downarrow}^{\ast} (Y)= \lim_{\rightarrow} \mathcal{E}_{B}^{\ast}$, where of course $B$ now ranges among the backward sets.

\begin{LE}\label{LE:duality-current-forms} The dual of the spaces of forward and compact/forward forms are
\[
\left(\mathcal{E}_{\uparrow}^{\ast}\left(  Y\right)\right)^{\prime}=\mathcal{D'}%
_{c \downarrow}^{\ast}\left(  Y\right)  \, \, and \, \, \left(\mathcal{E}_{c \uparrow}^{\ast}\left(  Y\right)\right)^{\prime}=\mathcal{D'}
_{\downarrow}^{\ast}\left(  Y\right) 
\]
i.e. the continuous functionals over the forward (resp. c/f) supported forms are exactly the compact/backward (resp. backward) supported currents. 
\end{LE} 

\noindent\textbf{Proof}. We will just consider the case of c/f support. Since the inclusion map $\mathcal{E}_{cpt}^{\ast
}  \hookrightarrow\mathcal{E}_{c\uparrow}^{\ast} $ is continuous with dense range, the adjoint map $\lp \mathcal{E}_{c \uparrow}^{\ast} \rp  ^{\prime}\rightarrow\mathcal{D}^{\prime}{}^{\ast}  $ is $1$-$1$ (and also weakly$*$ continuous, with weakly$*$ dense range), cf.[KN], p. 204. In particular each continuous linear functional on
$\mathcal{E}_{c\uparrow}^{\ast}  $ is a current, i.e. $\lp \mathcal{E}_{c \uparrow}^{\ast} \rp  ^{\prime} \subset \mathcal{D}^{\prime}{}^{\ast}$. 

If
$R\in\mathcal{D}_{\downarrow}^{\prime}{}^{\ast} $ has backward
support, then $R\left(  \varphi\right)  $ is defined for all $\varphi
\in\mathcal{E}_{c\uparrow}^{\ast} $, since spt$(R)\cap
$spt$\left(  \varphi\right)  $ is compact. Also, $R:\mathcal{E}_{A}^{\ast
}  \longrightarrow\mathbb{R}$ is continuous for any c/f set $A$. That is,
$\mathcal{D}_{\downarrow}^{\prime\ast}  \subset \lp \mathcal{E}%
_{c\uparrow}^{\ast}\rp ^{\prime}$.

Suppose instead $R\in \lp \mathcal{E}_{c\uparrow}^{\ast}\rp ^{\prime
}\subset\mathcal{D}^{\prime\ast}  $. If spt$\left(  R\right)
\cap A$ is not compact for some c/f supported $A$, then there exists a sequence of points in $A$ converging to $\infty$. We can thus find a larger c/f set $A'$, a sequence of pairwise disjoint balls $U_{n} \subset A'$ and a sequence of forms $\varphi_{n}\in\mathcal{E}_{cpt}^{\ast}(U_n) $ with $R\left(  \varphi_{n}\right)  =1$. Since $\varphi=\sum \varphi_{n}\in\mathcal{E}_{A'}^{\ast}  $ but $R\left(\varphi\right) $ is not defined, this is a contradiction. Therefore spt$\left(  R\right)  \cap A$ is compact, proving that $R \in \mathcal{D}_{\downarrow}^{\prime\ast} $ $\Box$\\

The previous lemma allows us to equip the spaces of currents $\mathcal{D'}_{c \uparrow} $ and $\mathcal{D'}_{\downarrow} $with the strong topology (dual to the space of forms with the right supports). Note that for any c/f set $A$, the subspace $\mathcal{D'}_{A} \subset \mathcal{D'}_{c \uparrow} $ is closed, as well as  $\mathcal{D'}_{B} \subset \mathcal{D'}_{\downarrow} $ for any backward set $B$. 

Finally, we endow the complexes $\mathcal{S}_{c \uparrow} \subset \mathcal{D'}_{c \uparrow} $ and $\mathcal{U}_{\downarrow} \subset \mathcal{D'}_{\downarrow}$ with the relative topology.

The complexes $\mathcal{S}_{c \uparrow}^{\ast}$ and $\mathcal{U}_{\downarrow}^{\ast}$ are fundamental here for several reasons: first, since they are made up of ``fewer'' objects, their topology is simpler; second, they are dual of each other; third, the $\mathcal{S}_{c \uparrow}$-complex is finitely dimensional over the Novikov ring $\Lambda$ (even if the $\mathcal{U}_{\downarrow}$-complex is not!).

\begin{LE}\label{LE:inductive-relative}
The space $\mathcal{S}_{c \uparrow}^{\ast}(Y)= \lim_{\rightarrow} \mathcal{S}_{A}^{\ast}$ (the inductive limit taken over the directed family of c/f sets $A$) and the topology described above equals the inductive limit topology. Similarly for $\mathcal{U}^{\ast}_{\downarrow} =  \lim_{\rightarrow} \mathcal{U}_{B}^{\ast}$.

Such topologies coincide respectively with the following two double limits topologies (the limits can be interchanged both times):
\begin{equation*}
\lim_{\overset{\longrightarrow}{m \rightarrow -\infty}} \lim_{\overset{\longleftarrow}{n \rightarrow +\infty}} \bigoplus_{m \leq f(y) \leq n } \Rs [ S_y] \,\,\, and \,\,\,
\lim_{\overset{\longrightarrow}{n\rightarrow +\infty}} \lim_{\overset{\longleftarrow}{m \rightarrow -\infty}} \prod_{m \leq f(y) \leq n } \Rs [ U_y] 
\end{equation*}

\end{LE}

\noi \textbf{Proof} We just consider the case of $\mathcal{S}_{c \uparrow}^{\ast}$. The currents defined by a single stable manifold $S_{y}$ form a discrete set in $\mathcal{D'}_{c \uparrow}^{\ast}$ (though they are not even a closed set in $\mathcal{D'}^{\ast}$ !) and hence they behave in $\mathcal{D'}_{c \uparrow}^{\ast}$ as the family of Dirac's point masses (or delta functions) supported at the critical points. 

The relative topology $\mathcal{S}_{c \uparrow}^{\ast} \subset \mathcal{D'}_{c \uparrow}^{\ast}$ is then represented by the double limit topology in the statement. On the other hand, if $A$ is a c/f set, the subspace $\mathcal{S}_{A}^{\ast}$ is topologically isomorphic to the product $\prod_{ S_x \subset A } \Rs [ S_x]$  and thus the inductive limit topology determined by the subspaces $\mathcal{S}_{A}^{\ast}$ also coincides with the double limit above $\Box$

\begin{LE}\label{LE:closed-open}
Any $\Lambda$-line in $\mathcal{S}_{c \uparrow}^{\ast}$ is a closed subspace. Moreover, $\mathcal{S}_{c \uparrow}^{\ast}$ is topologically isomorphic to a finite direct sum of such lines; in particular any $\Lambda$-linear subspace is closed and any $\Lambda$-linear self map of $\mathcal{S}_{c \uparrow}^{\ast}$ is relatively open, with closed range.
\end{LE}

Note that a $\Lambda$-line in $\mathcal{U}_{\downarrow}^{\ast}$ is not a closed subspace.\\

\noi \textbf{Proof}. Since multiplication by a nonzero $\lambda \in \Lambda$ is a topological isomorphism ($\Lambda$ being a field and multiplication being continuous), one can reduce to the case of the $\Lambda$-line generated by some $S_{y_0}$, so that the first statement is trivial. Now choose a (finite) $\Lambda$-basis for $\mathcal{S}_{c \uparrow}^{\ast}$, made up of elements like $S_{y_0}$: since the projections onto these $\Lambda$-lines are clearly continuous, the space $\mathcal{S}_{c \uparrow}^{\ast}$ is topologically isomorphic to a finite direct sum of such lines and the other statements follow $\Box$

\begin{PROP}\label{PROP:duality-U-S-complexes} The spaces $\mathcal{S}_{c \uparrow}^{\ast}(Y)$ and $\mathcal{U}_{\downarrow}^{\ast}(Y)$ are strong dual of each other under the current intersection pairing. The operator on $\mathcal{U}_{\downarrow}^{\ast}(Y)$, given by the adjoint of the differential $d$ on $\mathcal{S}_{c \uparrow}^{\ast}(Y)$, is the current boundary $\partial$. 
\end{PROP}

\noi \textbf{Proof}. Using Lemma \ref{LE:closed-open}, one can reduce to study the case when $\mathcal{S}_{c \uparrow}^{\ast}$ is 1-dimensional over $\Lambda$, i.e. to prove the duality between the space of c/f supported formal series in $q$ variables and the backward series in $q$ variables (forward and backward are supposed defined by a linear form in the lattice $\Zs^q$). But this case is conceptually the same as the duality between polynomials and formal series and we refer to [Tr], page 227 for a proof.

As for the adjoint of $d$, by smoothing, one easily checks that
\[
 U \bullet dS = \lim_{n \rightarrow +\infty} U \bullet d\beta_n = \lim_{n \rightarrow + \infty} \partial U \bullet \beta_n = \partial U \bullet S \,\,\,\,\, \Box
\]

Lemma \ref{LE:closed-open} proves that the range $d(\mathcal{S}_{c \uparrow}^{\ast})$ is a closed subspace and, combined with Proposition \ref{PROP:duality-U-S-complexes}, it also proves that the range $d(\mathcal{U}_{\downarrow}^{\ast})$ is closed, since the adjoint of a relatively open map has (weakly$*$ closed and hence) closed range, cf. [KN], p. 204,5. \\

The closed range property for the differentials $d$ in the $\mathcal{S}$ and $\mathcal{U}$-complexes is equivalent to the corresponding cohomology spaces $H^k(\mathcal{S}_{c \uparrow}^{\ast})$ and $H^k(\mathcal{U}_{\downarrow}^{\ast})$ being Hausdorff spaces (when endowed with the quotient topology).

Even more important, the same closed range property proves the nondegeneracy of the pairing in the following result, which is a standard corollary of Proposition \ref{PROP:duality-U-S-complexes} and the Hahn-Banach theorem:

\begin{PROP}\label{PROP:duality-U-S}
The current intersection between $\mathcal{S}_{c \uparrow}^{\ast}(Y)$ and $\mathcal{U}_{\downarrow}^{\ast}(Y)$ induces a pairing in cohomology:
\[
H^k( \mathcal{S}_{c \uparrow}^{\ast} (Y) ) \times H^{n-k} (\mathcal{U}_{\downarrow}^{\ast} (Y)) \rightarrow \Rs
\]
under which those locally convex Hausdorff spaces are in duality. 
\end{PROP}

\bigskip
\noi \textbf{Proof of Theorem \ref{THE:top.isomorphism}} We start by proving that the isomorphisms in formula (\ref{eq:THE:top.isomorphism}) are topological isomorphisms. Consider the operators, introduced in section \ref{SEC:Forward}), $\mathcal{E}^{\ast}_{c \uparrow} \overset {\mathbf{P}}{\rightarrow}  \mathcal{S}^{\ast}_{c \uparrow} \overset {\mathbf{I}}{ \hookrightarrow}  \mathcal{D'}^{\ast}_{c \uparrow}$, which induce the isomorphisms in cohomology. 
These operators are continuous (for the operator $\mathbf{P}$ cf. formula (\ref{eq:P})). Therefore, the maps induced in cohomology are continuous. If we show that these operators are open maps, the induced maps will be open as well, and we would be done. (Note that the inclusion $\mathcal{E}^{\ast}_{c \uparrow} \hookrightarrow \mathcal{D'}^{\ast}_{c \uparrow}$, is not an open map!).\\ 

The inclusion $\mathcal{S}^{\ast}_{c \uparrow} \hookrightarrow \mathcal{D'}^{\ast}_{c \uparrow}$ is open onto its image by definition (the topology on $\mathcal{S}^{\ast}_{c \uparrow}$ is the relative topology). 

Using formula (\ref{eq:P}), it is simple to prove that $\mathbf{P}: \mathcal{E}^{\ast}_{c \uparrow} \rightarrow  \mathcal{S}^{\ast}_{c \uparrow}$ is onto. Actually, $\mathbf{P}( \mathcal{E}_{A}^{\ast}) =\mathcal{S}_{A}^{\ast}$ for any c/f set $A\subset Y$ such that $\phi_t(A)\subset A$ for all $t\geq 0$, and no critical point belongs to $\partial A$ (such sets are cofinal in the c/f family). The open mapping theorem for Frechet spaces implies that the onto maps $\mathbf{P} : \mathcal{E}_{A}^{\ast} \rightarrow \mathcal{S}_{A}^{\ast}$ are open. If $U$ is a neighborhood of $0$ in $\mathcal{E}^{\ast}_{c \uparrow}$, then $ U \cap \mathcal{E}_{A}^{\ast}$ is a neighborhood of $0$ in $\mathcal{E}_{A}^{\ast}$ and hence $\mathbf{P}( U \cap \mathcal{E}_{A}^{\ast})$ (as well as $\mathbf{P}( U) \cap \mathcal{S}_{A}^{\ast} $, a fortiori) is a neighborhood of $0$ in $\mathcal{S}_{A}^{\ast}$. It follows that $\mathbf{P}( U)$ is a neighborhood of $0$ in $\mathcal{S}^{\ast}_{c \uparrow}$, i.e. $\mathbf{P}$ is open.
The proof that the isomorphisms in formula (\ref{eq:THE:top.isomorphism2}) are homeomorphisms is similar.\\ 

It remains to prove the duality statement. Since we proved that the spaces $H^k( \mathcal{S}_{c \uparrow}^{\ast} (Y) )$ and $H^{n-k} (\mathcal{U}_{\downarrow}^{\ast} (Y)) $ are Hausdorff, and the isomorphisms in formulae (\ref{eq:THE:top.isomorphism}) and (\ref{eq:THE:top.isomorphism2}) are homeomorphisms, it follows that the cohomology spaces defined by forms and currents are Hausdorff as well. This means that the differential $d:\mathcal{E}^{\ast}_{c \uparrow} \rightarrow \mathcal{E}^{\ast}_{c \uparrow} $ and the other differentials on the complexes of forms and currents are operators with closed ranges. Then, since the complexes of forms and currents (with the right supports) are dual, the duality on the corresponding cohomology spaces follows by the same Hahn-Banach argument invoked for Proposition \ref{PROP:duality-U-S}. The proof is complete $\Box$\\
 
We end by remarking that the geometric framework of Novikov theory was fundamental here. For a general Morse function $f$ on a manifold $Y$, one needs extra hypotheses to replace the arguments using $\Lambda$-finite dimensions. See [HH] for an example of what can go wrong with regards to duality for general differential operators.\\

\section*{Appendix}

We here collect and prove some algebraic properties of the Novikov ring and other rings of power series, needed in the present work. They were discovered by J.C. Sikorav and A. Pajitnov; for more informations we refer to Pajitnov's paper [P2]. 

We consider at once the general case of the Novikov ring $\Lambda$ in $q$ variables defined by the injective homomorphism $\chi : \Zs^q \rightarrow \Rs$ (see Definition \ref{DEF:Novikov-Ring}). The notations $\Lambda_{\Zs}$ or $\Lambda_{\Rs}$ will be used to specify the coefficients.

\bigskip

\noindent\textbf{Notations}. The \textbf{degree} of a monomial term $t^{n}$ is
defined to be $\chi \left(  n\right) \in \Rs$. The set of degrees of all the non zero monomial terms in the expansion of the series $\alpha\in\Lambda$ is denoted by $DEGS \left(  \alpha\right) \overset{def}{=} \chi\left(  \text{spt} (\alpha)   \right) \subset \Rs$. 

The set $DEGS\left(  \alpha\right) $ is a discrete (possibly finite) subset of
$\mathbb{R}$ , bounded below. Consequently, each $\alpha\in\Lambda$ has a
unique expansion $\alpha=\sum\limits_{j=0}^{n}a _{j}t^{A_{j}}$ with $n\in \Ns \cup \infty$, each $a_{j} \neq 0$ and $\deg t^{A_{j}}<\deg t^{A_{j+1}}$. The
\textbf{degree} of $\alpha$ is defined to be the degree of the \textbf{leading monomial} $a_{0}t^{A_{0}}$. The map $l:\Lambda\rightarrow\mathbb{Z}$ defined by
taking the \textbf{leading coefficient} $l\left(  \alpha\right)  =a_{0}$ is not a ring homomorphism; in fact $l(\alpha) + l(\beta)$ vanishes if the leading terms cancel. However $l(\alpha) l(\beta)$ always equals $l(\alpha \beta)$, and if $I$ is an ideal in $\Lambda$, then $l(I)$ is an ideal in $\mathbb{Z}$. 

\bigskip

\noindent\textbf{Lemma }\emph{An element }$\alpha\in\Lambda$\emph{ is a unit
if and only if }$l\left(  \alpha\right)  $\emph{ is a unit.}

\bigskip

\noindent\textbf{Proof.} We can assume $\alpha \!= \! 1 \! - \! \beta$ with $\deg\left(
\beta\right) \! > \! 0$. Since $\deg \left(  \beta^{k}\right)  \! = \! k\deg\left(
\beta\right)  $, the geometric series $1 + \beta + \beta^2 \ldots $ provides the inverse for
$\alpha$ in $\Lambda$   $\Box$\\

As an immediate corollary, we have:

\bigskip

\noindent\textbf{Algebraic Fact 1}\emph{ The Novikov ring with real coefficients }$\Lambda_{\mathbb{R}}
$\emph{ is a
field}.

\bigskip

The proof of the second algebraic fact is more involved:

\bigskip

\noindent\textbf{Algebraic Fact 2} \emph{The Novikov ring} $\Lambda_{\mathbb{Z}}$ \emph{is a principal ideal domain.}

\bigskip

\noindent\textbf{Proof}. Suppose $I$ is an ideal of $\Lambda_{\Zs}$. Since
$\mathbb{Z}$ is a p.i.d., $l\left(  I\right)  =\mathbb{Z}a$ for some integer
$a\in\mathbb{Z}$. Choose an element $\alpha=a+\sum\limits_{j=0}^{\infty}%
a_{j}t^{A_{j}}$ in the ideal $I$ with degree zero and leading coefficient
$l\left(  \alpha\right)  =a$. Given $\gamma\in I$, we will inductively define
$\beta=\sum\limits_{j=0}^{\infty}\beta_{j}\in\Lambda_{\mathbb{Z}}
$ so that $\gamma=\beta\alpha$, proving that $I=\Lambda_{\Zs} \alpha$.

Define $\gamma_{0}=\gamma$ and, given $\gamma_{k}\in I$, define the monomial
$\beta_{k}=b_{k}t^{B_{k}}$ as the leading term of $\gamma_{k}$ divided by $a$.
Since $l\left(  I\right)  =\mathbb{Z}a$, the coefficient $b_{k}\in\mathbb{Z}$.
Now define%
\[
\gamma_{k+1}=\gamma_{k}-\beta_{k}\alpha=\gamma-\left(  \beta_{0}+\beta
_{1}+..+\beta_{k}\right)  \alpha
\]
as the error in the factorization. Thus $\deg\beta_{k}=\deg\gamma_{k}%
<\deg\gamma_{k+1}$. Put $z_{k}=\deg\beta_{k}$: if no $\gamma_{k}$ vanishes, it
remains to show that $\lim\limits_{k\rightarrow\infty}z_{k}=+\infty$. Note
that $z_{k}=\min DEGS\left(  \gamma_{k}\right)  $, the set of degrees of terms
in $\gamma_{k}$. Let
\begin{align*}
\left\{  y_{1},y_{2},...\right\}   &  =\left\{  \deg t^{A_{1}},\deg t^{A_{2}%
},...\right\}  =DEGS\left(  \alpha\right)  \backslash\left\{  0\right\} \\
\left\{  x_{0},x_{1},...\right\}   &  =\left\{  \deg t^{C_{0}},\deg t^{C_{1}%
},...\right\}  =DEGS\left(  \gamma\right)  =DEGS\left(  \gamma_{0}\right)
\end{align*}
Now,

\vspace{-6mm}
\[
\gamma_{1}=\gamma_{0}-\beta_{0}\alpha=C_{1}t^{C_{1}}+C_{2}t^{C_{2}}%
+...-b_{0}a_{1}t^{C_{0}+A_{1}}-b_{0}a_{2}t^{C_{0}+A_{2}}-...
\]
and $z_{0}=x_{0}=\deg t^{C_{0}}$. Therefore,
\vskip-2mm
\[
DEGS\left(  \gamma_{1}\right)  \subset\left(  DEGS\left(  \gamma_{0}\right)
- \{ z_{0} \} \right)  \cup\left\{  z_{0}+y_{1},z_{0}+y_{2},...\right\}
\]
Similarly,
\vspace{-2mm}
\[
DEGS\left(  \gamma_{k+1}\right)  \subset\left(  DEGS\left(  \gamma_{k}\right)
- \{ z_{k} \} \right)  \cup\left\{  z_{k}+y_{1},z_{k}+y_{2},...\right\}
\]
\vskip-1mm
Consequently, the union of all the sets $DEGS\left(  \gamma_{k}\right)  $ is
contained in the set $D$ of real numbers of the form $x_{i_0}+y_{i_{1}%
}+...+y_{i_{k}}$. Since both the set of $x_{j}$'s and the set of $y_{i}$'s are
discrete and bounded below, the set $D$ is also discrete and bounded below.
Therefore $\lim\limits_{k\rightarrow\infty}z_{k}=\infty$ $\Box$

\bigskip

\noi \textbf{Remark: The Novikov ring $\Lambda_{\Zs}$ is Euclidean.} The norm is the absolute value of the leading coefficient, by the previous proof. 

\bigskip

We are left with the third algebraic fact:

\bigskip

\noindent\textbf{Algebraic Fact 3} \emph{The Novikov ring $\Lambda_{\mathbb{Z}}$ is flat over the Laurent Polynomial ring }$L=\Zs \left[  t_{1}
,{}_{\ldots},t_{q},t_{1}^{-1},{}_{\ldots},t_{q}^{-1}\right]$.

\bigskip

For the proof we need to introduce a new ring, also defined by Novikov in [N2]. \\

\noindent\textbf{Definition } \emph{A subset }$F\subset\mathbb{Z}^{q}$\emph{
in the lattice is \textbf{conical} (with respect to }$\chi$\emph{) if there exist }$a\in\mathbb{R}$\emph{ and }$\varepsilon>0$ \emph{s.t. }$F\subset\chi
^{-1}([a,+\infty))$\emph{ and (''stability'') this remains true for
all }$\chi^{\ast} \in \Rs^q$\emph{ with} $\left|  \chi-\chi^{\ast}\right|  <\varepsilon$.

\vspace{1ex}

\emph{The \textbf{Novikov conical} ring} ${}_{_<}\! \Lambda$\emph{ consists of all formal Laurent series }$\alpha
=\! \sum\limits_{n\in F} \! a_{n}t^{n}$\emph{ with integer or real coefficients whose support
}$F\! =$spt$(\alpha)  $\emph{ is a conical set in the lattice
}$\mathbb{Z}^{q}$.

\bigskip

Note that any conical set is compact/forward, so that the Novikov conical ring
${}_{_<}\!\Lambda$ is a subring of $\Lambda$. Again, the geometric
power series argument shows that the Novikov conical ring over the reals, ${}_{_<}\!\Lambda_{\Rs}$, is a field. \\

\noi \textbf{Proposition} \emph{ The Novikov conical ring over $\mathbb{Z}$}, ${}_{_<}\!\Lambda_{\mathbb{Z}} $\emph{, is a p.i.d.}\\

\noindent\textbf{Proof}. Given an ideal $I\subset{}_{_<}\!\Lambda\left(
\mathbb{Z}\right)  $, let $\bar{I}$ be the ideal in $\Lambda_{
\mathbb{Z}}  $ generated by $I$. Pick $\alpha\in I$, such that the ideal
$l\left(  I\right)  \subset\mathbb{Z}$ is generated by $l\left(
\alpha\right)  $. Then $\bar{I}=\Lambda_{\Zs} \alpha
$. In particular, if $\gamma\in I$ then $\beta=\gamma\alpha^{-1}\in
\Lambda_{\Zs} \subset \Lambda_{\Rs} $. Moreover, $\gamma \alpha^{-1} \in {}_{_<}\!\Lambda_{  \mathbb{R} } $ and since ${}_{_<}\!\Lambda
_{\Zs} ={}_{_<}\!\Lambda_{ \mathbb{R}}
\cap \Lambda_{\Zs} $, it follows $\beta \in {}_{_<}\!\Lambda_{\mathbb{Z}}$ $\Box$
\vskip4mm

\noindent\textbf{Theorem} \emph{The Novikov conical ring }${}_{_<}\!\Lambda_{\mathbb{Z}}$\emph{ is
a flat algebra over the Laurent Polynomial ring }$L=\Zs \left[  t_{1}
,{}_{\ldots},t_{q},t_{1}^{-1},{}_{\ldots},t_{q}^{-1}\right]$. 
\vskip4mm

\noindent\textbf{Proof}. Let ${}_{_<}\!\Lambda^{+}$ be the subring of the Novikov conical ring ${}_{_<}\!\Lambda$ made up of the series of positive degree, and $L^{+}$
the subring of $L$ made up of polynomials of positive degree (according to our
definition of degree), so that ${}_{_<}\!\Lambda^{+}$ is an algebra over $L^{+}$. The ring ${}_{_<}\!\Lambda$ can be presented as ${}_{_<}\!\Lambda = L\otimes_{L^{+}}{}_{_<}\!\Lambda^{+}$. Since "extending the scalars" preserves
flatness (cf. [B], I.2.7, Corollary 2), it's then enough to show that ${}_{_<}\!\Lambda^{+}$ is flat over $L^{+}$. 

Choose any $q$-ple $e=(e_1, .., e_q)$ with each $e_k \in
\Zs^q$ and $t^{e_k} \in L^+$ (i.e. $\chi (e_k)\! > \!0$), and introduce the following subsets of $\Rs^q$ and of the lattice $\Zs^q$ :

\begin{eqnarray*}
C(e)= \{ \sum x^k e_k \in  \Rs^q \, | \, x^k \in \Rs^+    \} &,& C^{1}(e)= \{ \sum x^k e_k  \in  \Rs^q \, | \, 0 \leq x^k \leq 1 \} \\
C_{\Zs}(e)= C(e) \cap \Zs^q  \, &,& \, 
C_{\Zs}^{1}(e)= C^{1}(e) \cap \Zs^q
\end{eqnarray*}

Note that $C_{\Zs}(e)$ might be bigger than the
integer span of the vectors in $e$, though it is always the
integer span of the vectors in $C_{\Zs}^{1}(e)$, which is a finite
set. Now define $\Lambda (e) \subset {}_{_<}\!\Lambda^+$ as the subring made up
of series with support in the cone $C(e)$ and $L(e) = L^+ \cap
\Lambda (e)$ the corresponding subring of polynomials. Equip the set of $q$-ples with the partial order defined by $e \leq e'$ if $C(e)
\subset C(e')$. It's evident that the
inductive limits ${}_{_<}\!\Lambda^+ =\lim\limits_{\longrightarrow} \Lambda(e)$ and $L^+ =\lim\limits_{\longrightarrow} L(e)$ hold.
Since the direct limit operation preserves flatness (cf. [B] I.2.7, Proposition 9 for the precise statement), the ring ${}_{_<}\!\Lambda^{+}$ is flat over $L^{+}$, provide we show that each ring of series $\Lambda(e)$ is flat over the corresponding ring of polynomials $L(e)$.

Fix $e$ as above. As already observed, $L(e)$ is generated by the monomials with support in $C_{\Zs}^{1}(e)$, which are a finite number, say $r$. The ring $L(e)$ is hence Noetherian since there exists some surjective homomorphism $\Zs [x_1 , .., x_r] \rightarrow L(e)$ (any quotient of a Noetherian ring is Noetherian).

Consider now the ideal $\mathcal{M} \subset L(e)$ consisting of the polynomials with vanishing free term. The ring $\Lambda(e)$ is the completion
of $L(e)$ for this ideal ($\Lambda(e) = \lim\limits_{\overset{\longleftarrow}{\infty \leftarrow k}} L(e) / \mathcal{M}^k$ ). Using [B] III, 3.4, Theorem 3, this proves that $\Lambda(e)$ is flat over $L(e)$, completing the proof  $\Box$\\

\noindent\textbf{Proof of Fact 3} Any torsion free module over a p.i.d. is flat ([B], I.2.4 Proposition 3), and the Novikov ring $\Lambda$ is torsion free over the Novikov conical
ring ${}_{_<}\!\Lambda$, so $\Lambda$ is flat over ${}_{_<}\!\Lambda$. We just
proved ${}_{_<}\!\Lambda$ is flat over the Laurent polynomials $L$ and since
flatness is a transitive property, the result follows $\Box$

\begin{center}
{\footnotesize REFERENCES}
\end{center}
\small{

\noindent \scriptsize{[AB]} \scriptsize{D. AUSTIN and P.J. BRAAM}, Morse Bott Theory and equivariant cohomology, in \emph{The Floer Memorial Volume}, Birkhauser 1994\newline
[B] \scriptsize{ N. BOURBAKI}, \emph{Commutative Algebra Chap.1-7} Springer-Verlag\newline
[Bo] \scriptsize{R. BOTT}, Morse Theory Indomitable, \emph{Publ. Math. IHES} \textbf{68} (1988) 99-114 \newline
[D]\scriptsize{ G. DE RHAM}, \emph{Differentiable manifolds. Forms,
currents, harmonic forms}. GMW \textbf{266}. Springer-Verlag 1984\newline
[F]\scriptsize{ M. FARBER}, Exactness of the Novikov inequalities, \emph{Functional Anal. Appl.} \textbf{19} (1985) 40-48\newline
[FR] \scriptsize{M. FARBER and A. RANICKI}, The Morse Novikov theory of circle valued functions and noncommutative localization, \emph{Tr. Math. Inst. Steklova} \textbf{225} (1999) 381-388\newline
[Fe] \scriptsize{H. FEDERER}, \emph{Geometric measure theory}, GMW \textbf{153}
Springer-Verlag 1969\newline
[Fl] \scriptsize{A. FLOER}, Witten's complex and infinite dimensional Morse theory, \emph{J. diff. geom.} \textbf{30} (1989) 207-221\newline  
[G] \scriptsize{R. GODEMENT}, \emph{Topologie alg\'{e}brique et th\'{e}orie des faisceaux.} Hermann, Paris 1973\newline
[HH] \scriptsize{C. HARVEY and F.R. HARVEY}, Open mappings and the lack of fully completeness of $\mathcal{D'}(\Omega )$, \emph{Proc. AMS} \textbf{25} (1970) 786-790\newline 
[HL] \scriptsize{F.R. HARVEY and H.B. LAWSON}, Finite volume flows and Morse theory, \emph{Ann.of Math}. \textbf{153} (2001) 1-25\newline
[HP] \scriptsize{F.R. HARVEY and J. POLKING}, Fundamental solutions in complex analysis, Part I, \emph{Duke Math J}. \textbf{46} (1979) 253-300\newline
[HZ] \scriptsize{F.R. Harvey and J. Zweck} Stiefel-Whitney currents, \emph{J.Geom.Anal.} \textbf{8} (1998), no.5 805-840 \newline
[KN] \scriptsize{J. KELLEY and I. NAMIOKA}, \emph{Linear Topological Spaces} Van Nostrand 1963 \newline
[L] \scriptsize{F. LATOUR}, Existence de 1-formes ferm\'ee non singuli\`eres dans une classe de cohomologie de de Rham, \emph{Publ. Math. IHES} \textbf{80} (1994) 135-194 \newline 
[La] \scriptsize{F. LAUDENBACH}, On the Thom-Smale complex, (appendix) Asterisque \textbf{205} (1992) \newline 
[Lt] \scriptsize{J. LATSCHEV}, Gradient flows of Morse-Bott functions,\emph{Math.Ann.} \textbf{318} (2000) no.4, 731--759\newline
[M] \scriptsize{G. MINERVINI}, PhD thesis, Universit\`a ``La Sapienza'', Roma 2003.\newline
[M2] \scriptsize{G. MINERVINI}, Preprint.\newline
[N] \scriptsize{S.P. NOVIKOV}, Multivalued functions and functionals. An analogue of the Morse theory, \emph{Sov. Math. Dolklady} \textbf{24} (1981) 222-226
\newline
[N2] \scriptsize{S.P. NOVIKOV}, 
The Hamiltonian formalism and a multivalued analogue of Morse theory
\emph{Russ. Math. Survey} \textbf{37} (1982) no.5 1-56\newline
[P] \scriptsize{A.V. PAJITNOV}, Counting closed orbits of Gradients of Circle
Valued Maps, arXiv preprint DG/0104273 v2 (2002)\newline
[P2] \scriptsize{A.V. PAJITNOV}, On the sharpness of Novikov type inequalities for manifolds with free abelian fundamental group. \emph{Math. USSR Sbornik} \textbf{68} (1991) no.2, 351-389\newline
[P3] \scriptsize{A.V. PAJITNOV}, On the Novikov complex for rational Morse forms, preprint Odense University \textbf{12} (1991), reprinted in \emph{Ann.Fac.Sci.Toul.} \textbf{4} (1995), no.2 297-338 \newline
[R] \scriptsize{A. RANICKI}, Circle Valued Morse Theory and Novikov Homology, e-print AT.0111317 (2001)\newline
[S] \scriptsize{L. SHILNIKOV et al}, \emph{Methods of non qualitative theory in nonlinear dynamics Part I} Nonlinear Science, World Scientific 1998\newline
[Sc] \scriptsize{D. SCHUTZ} Gradient flows of closed 1-forms and their closed orbits, \emph{Forum Math.} \textbf{14} (2002) 509-703\newline
[Sc2] \scriptsize{D. SCHUTZ} Controlled connectivity of closed 1-forms, \emph{Alg. Geom. Top.} \textbf{2} (2002) 177-217\newline
[Sm] \scriptsize{S. SMALE}, On gradient dynamical systems. \emph{Ann.of Math}. 2, \textbf{74} (1961) 199-206;\newline
[T] \scriptsize{R. THOM} Sur une partition en cellules associ\'ee \`a une function sur une vari\'et\'e \emph{C.R. Acad. Sci. Paris} \textbf{228} (1949) 973-975 \newline[Tr] \scriptsize{F. TREVES}, \emph{Topological Vector Spaces, Distributions and Kernels}, Academic Press, 1967 \newline
[W] \scriptsize{E. WITTEN}, Supersymmetry and Morse Theory, \emph{J. diff. geom.} \textbf{17} (1982) 661-692\newline
\end{document}